%

\magnification=\magstep1
\def\forces{\parallel\!\!\! -}
\def\restrict{\lceil}
\def\Smallskip{\vskip1.5truecm}
\def\Bigskip{\vskip2.5truecm}

\def\qed{{\vcenter{\hrule height.4pt \hbox{\vrule width.4pt height5pt
 \kern5pt \vrule width.4pt} \hrule height.4pt}}}
\def\ok{\vbox{\hrule height 8pt width 8pt depth -7.4pt
    \hbox{\vrule width 0.6pt height 7.4pt \kern 7.4pt \vrule width 0.6pt height 7.4pt}
    \hrule height 0.6pt width 8pt}}
\def\nt{{\leq}\kern-1.5pt \vrule height 6.5pt width.8pt depth-0.5pt \kern 1pt}
\def\sd{{\times}\kern-2pt \vrule height 5pt width.6pt depth0pt \kern1pt}
\def\notin{{\in}\kern-5.5pt / \kern1pt}
\def\ZZ{{\rm Z}\kern-3.8pt {\rm Z} \kern2pt}
\def\RR{{\rm I\kern-1.6pt {\rm R}}}
\def\NN{{\rm I\kern-1.6pt {\rm N}}}
\def\QQ{{\rm Q}\kern-4.4pt {\rm \vrule height6.8pt 
    width.3pt depth-1pt} \kern7.4pt}
\def\CC{{\rm \kern 3.5pt \vrule height 6.5pt 
    width.3pt depth-1.3pt \kern-4pt C \kern.8pt}}
\def\PP{{\rm I\kern-1.6pt {\rm P}}}
\def\BB{{\rm I\kern-1.6pt {\rm B}}}
\def\KK{{\rm I\kern-1.6pt {\rm K}}}
\def\11{{\rm 1}\kern-2.2pt {\rm \vrule height6.1pt
    width.3pt depth0pt} \kern5.5pt}
\def\zp#1{{\hochss Y}\kern-3pt$_{#1}$\kern-1pt}

\def\egs{\vrule height 6pt width.5pt depth 2.5pt \kern1pt}
\font\capit=cmcsc10 scaled\magstep0
\font\capitg=cmcsc10 scaled\magstep1
\font\dunhg=cmcsc10 scaled\magstep1
\font\dunhgg=cmcsc10 scaled\magstep1
\font\sanseg=cmss10 scaled\magstep0
\font\sanse=cmss10 scaled\magstep0
\overfullrule=0pt
\openup1.5\jot

\centerline{\dunhgg Set-theoretic aspects of periodic FC-groups ---}
\Smallskip
\centerline{\dunhgg Extraspecial p-groups and Kurepa trees\footnote
{$^*$}{\rm This work forms part of the author's PhDThesis written under
the supervision of Prof. Ulrich Felgner, T\"ubingen, 1991.}}
\Bigskip
\centerline{\capitg J\"org Brendle}
\Smallskip
\centerline{\sanseg Mathematisches Institut} 
\centerline{\sanseg der Universit\"at T\"ubingen}
\centerline{\sanseg Auf der Morgenstelle 10}
\centerline{\sanseg W-7400 T\"ubingen}
\centerline{\sanseg Germany}
\Bigskip
\centerline{\capit Abstract}
\bigskip
\noindent When generalizing a characterization of centre-by-finite
groups due to B. H. Neumann, M. J. Tomkinson asked the following question.
Is there an $FC$-group $G$ with $\vert G / Z(G) \vert = \kappa$ but
$[ G : N_G(U) ] < \kappa$ for all (abelian) subgroups $U$
of $G$, where $\kappa$ is an uncountable cardinal [16, Question
7A on p. 149]. We consider this
question for $\kappa = \omega_1$ and $\kappa = \omega_2$. It turns
out that the answer is largely independent of $ZFC$ (the usual axioms
of set theory), and that it differs greatly in the two cases.
\vfill\eject
\noindent{\dunhg Introduction}
\Smallskip
{\sanse The problem and its history.} 
The purpose of this paper is to give some independence proofs concerning
the existence of $FC$-groups having some additional properties.
Recall that a group $G$ is $FC$ iff every element $g \in G$
has finitely many conjugates; i.e. iff $[G:C_G(g)]$ is finite
for any $g \in G$. We shall mostly be concerned with periodic 
$FC$-groups. \par
In the fifties, B. H. Neumann gave the following characterization
of {\it centre-by-finite groups}; i.e. groups with $\vert G / Z(G) \vert
< \omega$. \smallskip
{\it \item{(I)} The following are equivalent for any group $G$. \par
\itemitem{(i)} $G$ is centre-by-finite. \par
\itemitem{(ii)} Each subgroup of $G$ has only finitely many conjugates;
i.e. $[G : N_G (U)] < \omega$ for all $U \leq G$. \par
\item{} If $G$ is an $FC$-group both are equivalent to \par
\itemitem{(iii)} $U/U_G$ is finite for all $U \leq G$.}
\smallskip
\noindent Here $U_G$ denotes the largest normal subgroup of $G$ contained 
in $U$, i.e. $U_G := \cap_{g \in G} \; g^{-1} U g$; it is called the
{\it core} of $U$ in $G$. It was indicated by Eremin that in both
(ii) and (iii) above it suffices to consider abelian subgroups
(cf [16, 7.12(a) and 7.20]; also note that a group satisfying
(iii) above is in general not $FC$ [16, p. 142]). \par
Following M. J. Tomkinson [15] (see also [4]), 
for an infinite cardinal $\kappa$,
let ${\cal Z}_\kappa$ denote the class of groups $G$ in which
$[ G : C_G(U) ] < \kappa$ whenever $U \leq G$ is generated by fewer than
$\kappa$ elements (for $\kappa > \omega$,
this is equivalent to saying that $[ G : C_G(U) ] < \kappa$ for
$U \leq G$ of size less than $\kappa$). Clearly ${\cal Z}_\omega$
is just the class of $FC$-groups. Generalizing Neumann's result Tomkinson
proved in [15] (see also [16, Theorem 7.20]).
\smallskip
{\it \item{(II)} Let $\kappa$ be an infinite cardinal. The following
are equivalent for any $FC$-group $G$ in ${\cal Z}_\kappa$. \par
\itemitem{(i)} $\vert G / Z(G) \vert < \kappa$. \par
\itemitem{(ii)} $[G : N_G (U) ] < \kappa$ for all $U\leq G$. \par
\itemitem{(iii)} $[G : N_G(A) ] < \kappa$ for all abelian $A \leq G$. \par
\itemitem{(iv)} $\vert U / U_G \vert < \kappa$ for all $U \leq G$. \par
\itemitem{(v)} $\vert A /A_G \vert < \kappa$ for all abelian $A \leq G$.}
\smallskip
\noindent It was later shown by Faber and Tomkinson in [4]
that the condition that $G$ is $FC$ can be dropped in (II). \par
On the other hand one might ask whether the condition that $G$
is ${\cal Z}_\kappa$ is necessary to prove the equivalence of
(i) through (v) in (II). Clearly the following implications always
hold. \smallskip
\centerline{$ (i) \Longrightarrow (ii) \Longrightarrow (iii)$} \par
\centerline{$ (i) \Longrightarrow (iv) \Longrightarrow (v)$} \smallskip
\noindent But what about the others? \par
One answer to this question was indicated by Tomkinson himself
[16, p. 149]. Recall that a $p$-group $E$ is called {\it extraspecial}
iff $\Phi (E) = E' = Z(E) \cong \ZZ_p$. (Here $\Phi (E)$ is the {\it
Frattini subgroup} of $E$; i.e. the intersection of all maximal
subgroups of $E$.) This implies that $E/E'$ is elementary abelian.
Let $E$ be an extraspecial $p$-group of size $\omega_1$
all of whose abelian subgroups are countable (the existence of such groups
was proved by S. Shelah and J. Stepr\=ans in [13] improving on earlier
work of A. Ehrenfeucht and V. Faber [16, Theorem 3.12]
who got the same result under
the additional assumption of the continuum hypothesis ($CH$)).
Note that if $U \leq E$ then $U$ is either normal and so $U=U_E$
or $U$ is abelian and so $\vert U / U_E \vert \leq \vert U \vert \leq
\omega$. So for $\kappa = \omega_1$ and $G = E$, (iv) and (v) in (II)
are true, whereas (i) is not. Also, if $A \leq E$ is maximal (abelian)
with respect to $A \cap E' = 1$, then $\langle E' , A \rangle =
C_E(A) = N_E(A)$. Hence $\vert A \vert = \vert N_E(A) \vert = \omega$
and $[E : N_E (A) ] = \omega_1$. Thus (ii) and (iii) do not hold either.
\par
Why is this so? -- To get (iv) and (v) but not (i) we used an extraspecial
$p$-group such that all (maximal) abelian subgroups are small.
Dually, to get (ii) and (iii) but not (i) we {\it should} use an
extraspecial $p$-group such that all maximal abelian subgroups are
large in the sense that their indices are small. It will be one of
our goals to discuss the existence of such groups (see Theorems D and
E below and $\S$ 5). Tomkinson proved already in [14]
that there are no such groups of size $\omega_1$ (this also follows from
our more general Theorem C).
\bigskip

{\sanse The main results.} For $\kappa = \omega_1$ our results are as 
follows.
\smallskip  
{\capit Theorem A.} {\it Under $CH$ there is an $FC$-group $G$ with $\vert
G/Z(G) \vert = \omega_1$ but $[G:N_G(A)] \leq \omega$ for all abelian
subgroups $A \leq G$.} \smallskip
{\capit Theorem B.} {\it It is consistent (assuming the consistency of
$ZFC$) that there is no $FC$-group $G$ with $\vert
G/Z(G) \vert = \omega_1$ but $[G:N_G(A)] \leq \omega$ for all abelian
subgroups $A \leq G$.} \smallskip
\noindent Theorems A and B show that the question whether (i) and (iii)
in (II) are equivalent for $\kappa = \omega_1$ is not decided by the
axioms of set theory alone. The example used to prove Theorem A has a countable
subgroup $U$ with $[G:N_G(U)] = \omega_1$. So it does not answer
the following
\smallskip
{\capit Question 1.} {\it Let $\kappa = \omega_1$. Are (i) and (ii) in
(II) equivalent for all $FC$-groups $G$?}
\smallskip
\noindent We conjecture that the answer is yes. Our reason for believing
this is the following partial result.
\smallskip
{\capit Theorem C.} {\it Let $\kappa = \omega_1$. Then (i) through (iii)
in (II) are equivalent for all finite-by-abelian groups $G$.}
\smallskip
\noindent Here, a group $G$ is called {\it finite-by-abelian} iff $G'$ is finite.
\par
Question 1 and Theorem C are very closely related to another problem
of Tomkinson [16, Question 3F on p. 60]. Following [14] (see also
[16, chapter 3]) let
${\cal Z}$ be the class of locally finite groups $G$
satisfying: for all cardinals $\kappa$ and all $H \leq G$ of size
$<\kappa$, $[G:C_G(H)] < \kappa$. So ${\cal Z}$ is the class of
periodic groups in the intersection of the ${\cal Z}_\kappa$.
And ${\cal Y}$ is the class of locally finite groups $G$ satisfying:
for all cardinals $\kappa$ and all $H \leq G$ of size $<\kappa$,
$[G:N_G(H)] < \kappa$. Clearly ${\cal Z} \subseteq {\cal Y}$.
Tomkinson asked whether there are ${\cal Y}$-groups
which are not in ${\cal Z}$. He proved in [14, Theorem D(i)] that
any extraspecial $p$-group in ${\cal Y}$ of size $\omega_1$ 
lies in ${\cal Z}$. We generalize this by showing
\smallskip
{\capit Theorem C'.} {\it ${\cal Y} = {\cal Z}$ for finite-by-abelian
groups of size $\omega_1$.}
\smallskip
{\capit Theorem B'.} {\it Assuming the consistency of $ZFC$ it is
consistent that ${\cal Y} = {\cal Z}$ for $FC$-groups of size $\omega_1$.}
\smallskip
\noindent The proofs of these results use the same ideas as the proofs of 
Theorem C and B, respectively, and
we hope that our argument can be generalized to give a positive answer to
\smallskip
{\capit Question 1'.} {\it Is ${\cal Y} = {\cal Z}$ for $FC$-groups of size
$\omega_1$?}
\smallskip
\noindent Note that a positive answer to Question 1' would give a positive
answer to Question 1 too. For suppose there is a counterexample $G$.
Then $\vert G/Z(G) \vert = \omega_1$ (and without loss we may
assume that $\vert G \vert = \omega_1$) but $[G:N_G(U)] \leq \omega$ for
all $U \leq G$. By Tomkinson's result (II), $G \notin {\cal Z}_{\omega_1}$,
so $G \notin {\cal Z}$, hence $G \notin {\cal Y}$; i.e. there is
a countable $U \leq G$ such that $[G:N_G(U)] = \omega_1$, a contradiction. 
-- The problem seems to be of group theoretical character,
and might involve a better understanding of countable periodic
$FC$-groups.
\par
For $\kappa = \omega_2$ the picture changes considerably. Recall
that an uncountable cardinal $\kappa$ is {\it (strongly) inaccessible} iff it
is regular (i.e., it is not the union of $< \kappa$ sets of size
$< \kappa$) and all cardinals $\lambda < \kappa$ satisfy $2^\lambda 
< \kappa$
(especially, $\kappa$ is a limit cardinal). The existence of inaccessible
cardinals cannot be proved in $ZFC$; in fact, something much stronger
is true. Let $I$ denote the sentence {\it there is an inaccessible
cardinal}. Then the consistency of $ZFC$ can be proved in the system
$ZFC+I$ (see [9, chapter IV, Theorem 6.6 and p. 145]). Hence, by
G\"odel's Incompleteness Theorem, the consistency of $ZFC+I$ cannot be
proved
from the consistency of $ZFC$ alone. -- A {\it weak Kurepa tree}
is a tree of height $\omega_1$ with $\omega_2$ uncountable branches such that
all levels have size $\leq \omega_1$. A {\it Kurepa tree} is a
weak Kurepa tree with countable levels (a more formal definition will
be given in $\S$ 1). The existence of Kurepa trees is consistent
(assuming the consistency of $ZFC$), 
and the non-existence of Kurepa trees is equiconsistent with the existence
of an inaccessible (see, again, our $\S$ 1 for details). 
\smallskip
{\capit Theorem D.} {\it Assume there is a Kurepa tree. Then there is an
extraspecial $p$-group of size $\omega_2$ such that $[G:A] \leq \omega_1$
for all maximal abelian subgroups $A \leq G$.}
\smallskip
\noindent On the other hand one can show (Theorem 5.4) that the existence 
of such a group implies the existence of a weak Kurepa tree. 
In fact, we can prove the following much stronger result.
\smallskip
{\capit Theorem E.} {\it Assuming the consistency of $ZFC + I$ it is
consistent that for both $\kappa = \omega_1$ and $\kappa = \omega_2$ and
any $FC$-group $G$, (i) through (iii) in (II) are equivalent.}
\smallskip
\noindent We thus get
\smallskip
{\capit Corollary.} {\it The following theories are equiconsistent.
\par
\item{(a)} $ZFC + I$. \par
\item{(b)} $ZFC \; +$ for $\kappa = \omega_1$ and $\kappa = \omega_2$ and
for any $FC$-group $G$: if $\vert G / Z(G) \vert = \kappa$ then there
is an abelian subgroup $A \leq G$ with $[G : N_G(A)] = \kappa$. \par
\item{(c)} $ZFC \; +$ any extraspecial $p$-group of size $\omega_2$
has an (abelian) subgroup with $[G:N_G(A)] = \omega_2$. \par}
\smallskip
\noindent "$(a) \Rightarrow (b)$" is Theorem E; "$(b) \Rightarrow (c)$"
is trivial; and "$(c) \Rightarrow (a)"$ follows from Theorem D using the
equiconsistency concerning the non-existence of Kurepa trees  
mentioned above. It should be pointed out that
we cannot prove the {\it equivalence} of (b) and (c); namely, it is
consistent (assuming again the consistency of $ZFC + I$) that (c) holds
but (b) does not (see Theorem 5.8). Still our Corollary shows 
again (cf [14] or [16, chapter 3, especially 3.15]) the
importance of the extraspecial $p$-groups in the class of periodic 
$FC$-groups.
\bigskip

{\sanse The organization of the paper.} Our results use mainly classical
(modern) set theory. For algebraists who might not be familiar with this
material, we give a short Introduction to this subject in $\S$ 1.
We hope that this makes our work more intelligible. The reader who has
seen forcing etc. before should skip the entire $\S$ 1.
\par
In the second section we show that a countable finite-by-abelian group
is generated by finitely many abelian subgroups (Theorem 2.2). We also
discuss what goes wrong when {\it countable} is dropped from
the assumption of the Theorem.
\par
In the third section we prove a result on automorphisms of countable 
periodic abelian groups which turns out to be crucial for our
arguments (Theorem 3.2); we will apply it in the proofs of Theorems
B, C, and E. If we could generalize this result to countable periodic
$FC$-groups, we would get a positive answer to Questions 1 and 1'.
Our original proof involved a fragment of Ulm's classification
theorem. Since then, M. J. Tomkinson has found a much shorter
and more elegant proof which we reproduce with his permission...
We close $\S$ 3 with the proof of Theorem C (and C').
\par
Section 4 is devoted to the proofs of Theorems A and B (and B'); i.e.
to the case $\kappa = \omega_1$. It turns out that the knowledge of
maximal abelian subgroups of countable periodic $FC$-groups is essential.
\par
In section 5 we deal with the case $\kappa = \omega_2$ and the
relationship between Kurepa trees and extraspecial $p$-groups;
we prove Theorems D and E. We think that those results are the most
interesting and most beautiful of our work.
\par
We close with some generalizations in $\S$ 6.
\par
Finally note that we get most of the {\it main results} mentioned in
the preceding subsection as corollaries to more technical theorems
and constructions, and we hope that the ideas involved in the latter
might be useful when dealing with other problems as well. They are
Theorems 2.2, 3.2 and 4.4 (with its elaboration in 5.7) and the
easy 5.4 (see also 5.7) -- and the constructions in 4.2 (modified in
4.6 and 4.7) and 5.3 (modified in 5.9).
\bigskip

{\sanse Group-theoretic notation and basic facts on $FC$-groups.} 
Our group-theoretic 
notation is standard. Good references are [12] for general group
theory, [5] and [6] for abelian groups (which will be written
additively), and [16] for $FC$-groups.
\smallskip
For completeness' sake we give our extension-theoretic notation.
Let $A$, $G$ be groups. If $G \leq Aut (A)$, we let $A \sd G$ denote
the semidirect product of $A$ and $G$. If $\tau: G^2 \to A$ is
a factor system (i.e. $\forall g \in G \; (\tau (g,1) = \tau (1,g) = 1)$ 
and $\forall f,g,h \in G \; ( \tau (fg,h) \tau (f,g) = \tau (f, gh)
\tau(g,h))$), we let $E(\tau)$ denote the corresponding extension
(where the operation of $G$ on $A$ is trivial). In the latter case,
group multiplication is given by the formula 
$$\forall (a,g), (b,h) \in E(\tau) \;\;\; (a,g)*(b,h) = (\tau(g,h)ab,gh).$$
More details can be found in [12, chapter 11].
\smallskip
We note that in all of our results (in particular, in Theorems B, C, and E)
it suffices to consider {\it periodic} $FC$-groups. The reason for
this is as follows. By a result of \v Cernikov [16, Theorem 1.7],
any $FC$-group can be embedded in a direct product of a periodic
$FC$-group and a torsion-free abelian group. Now suppose $G$
is an (arbitrary) counterexample to one of our results; i.e.
$\vert G/Z(G) \vert = \kappa$, but $[G:N_G(A)] < \kappa$ for
all (abelian) $A \leq G$. Assume $G \leq P \times T$ and
$\pi (G) = P$, $\rho (G) =T$, where $P$ is a periodic $FC$-group
and $T$ is torsion-free abelian, and $\pi$ and $\rho$ are the projections.
Clearly $\vert P/Z(P) \vert = \kappa$, and also $[P:N_P(A)] < \kappa$
for all (abelian) $A \leq P$. This gives us a periodic counterexample.
\bigskip
{\sanse Acknowledgments.} I should like to thank Ulrich Felgner and
Frieder Haug for many stimulating discussions relating to the
material of this work. I am also grateful to M. J. Tomkinson
for simplifying the proof of Theorem 3.2, and to both him
and the referee for many valuable suggestions.
\Bigskip
\noindent{\dunhg $\S$ 1. Set-theoretic preliminaries}
\Smallskip
{\sanse Set-theoretic Notation.} If $X$ is a set, $[X]^\kappa$
denotes the set of subsets of $X$ of size $\kappa$; $[X]^{< \kappa}$
is the set of subsets of $X$ of size $< \kappa$; $[X]^{\leq \kappa}$
etc. are defined similarly. If $X \in [\kappa]^n$ for some $n \in \omega$,
then $X(i)$ $(i < n)$ denotes the $i$-th element of $X$ under the
inherited ordering. Further set-theoretic notation can be found in [9]
or [7].
\bigskip
{\sanse Delta-systems and almost disjoint sets.} A family ${\cal A}$
of sets is called a {\it delta-system ($\Delta$-system)} if there
is an $R$ (called the {\it root} of ${\cal A}$) such that
$$\forall A, B \in {\cal A} \; \; ( A \ne B \Longrightarrow A \cap
B = R ).$$
The {\it delta-system lemma} [9, chapter II, Theorem 1.6]  asserts
that given a collection ${\cal A}$ of sets of size $< \kappa$ with
$\vert {\cal A} \vert \geq \theta$ where $\theta > \kappa$ is regular
and satisfies $\forall \alpha < \theta \; ( \vert \alpha^{<\kappa} \vert
< \theta )$, there is a ${\cal B} \subseteq {\cal A}$ of size
$\theta$ which forms a $\Delta$-system. We shall use it most often in
case $\kappa = \omega$. \par
If $\kappa$ is a cardinal, a family of sets ${\cal A} \subseteq {\cal P}
(\kappa)$ is called {\it almost disjoint (a.d.)} iff
$$\forall A \in {\cal A} \; (\vert A \vert = \kappa) \; {\rm and} \;
\forall A, B \in {\cal A} \; (A \ne B \Rightarrow \vert A \cap B \vert
< \kappa).$$
\bigskip
{\sanse Trees.} A {\it tree} is a partial order $\langle T, \leq \rangle$
such that for each $x \in T$, $\{ y \in T ; \; y < x\}$ is wellordered
by $<$. Let $T$ be a tree. For $x \in T$, the {\it height} of $x$ in
$T$ ($ht(x,T)$) is the order type of $\{ y \in T ; \; y < x \}$.
For each ordinal $\alpha$, the {\it $\alpha$-th level} of $T$ is
$Lev_\alpha(T) = \{x \in T; \; ht(x,T)= \alpha \}$. The {\it height}
of $T$ ($ht(T)$) is the least ordinal $\alpha$ such that $Lev_\alpha (T)
= \emptyset$. A {\it branch} of $T$ is a maximal totally ordered
subset of $T$. \par
A {\it weak Kurepa tree} is a tree $T$ of height $\omega_1$ with at least
$\omega_2$ uncountable branches such that $\forall \alpha < \omega_1 \; ( \vert
Lev_\alpha (T) \vert \leq \omega_1 )$. Clearly, if $CH$ holds,
the complete binary tree of height $\omega_1$ is a weak Kurepa tree.
A {\it Kurepa tree} is a weak Kurepa tree $T$ satisfying $\forall \alpha
< \omega_1 \; ( \vert Lev_\alpha (T) \vert \leq \omega)$. A {\it
Kurepa family} is an ${\cal F} \subseteq {\cal P} (\omega_1)$ such that
$\vert {\cal F}\vert \geq \omega_2$ and $\forall \alpha < \omega_1
\; ( \vert \{ A \cap \alpha ; \; A \in {\cal F} \} \vert \leq \omega)$.
It is easy to see [9, chapter II, Theorem 5.18] that there is a 
Kurepa family iff there is a Kurepa tree.
\bigskip
{\sanse Partial orders and forcing.} Forcing was created by Cohen
in the early sixties to solve Cantor's famous continuum problem; i.e.
to show that for any cardinal $\kappa$ of cofinality $> \omega$ it
is consistent that $2^\omega = \kappa$ -- assuming the consistency of $ZFC$.
Since then many other independence problems have been solved by the same method.
As forcing will occupy a central position in our work, we briefly
define its main notions. For a
(very nicely written) introduction to this subject, we refer the
reader to [9]. \par
Let $\langle \PP , \leq \rangle $ be a partial order (p.o. for short;
sometimes, $\PP$ will be referred to as {\it forcing notion}).  
The elements of $\PP$ are called {\it conditions}. If $p, q \in \PP$
and $p \leq q$,
then $p$ is {\it stronger than} $q$ (or $p$ is said to {\it extend}
$q$). $p$ and $q$ are {\it compatible} iff $\exists r \in \PP
\; (r \leq p \; \land r \leq q)$; otherwise they are {\it incompatible}
($p \bot q$).
A set $D \subseteq \PP$ is called {\it dense} iff $\forall p \in \PP
\; \exists q \leq p \; (q \in D)$; $D$ is {\it open dense} iff it is
dense and $\forall p \in \PP \; \forall q \in D \; (p \leq q \Rightarrow
p \in D)$. ${\cal G} \subseteq \PP$ is called a {\it filter} iff
$\forall p, q \in {\cal G} \; \exists r \in {\cal G} \; ( r \leq p \;
\land \; r \leq q)$ and $\forall p \in {\cal G}\; \forall q \in \PP \;
(p \leq q \Rightarrow q \in {\cal G})$. -- Now suppose ${\cal M}$ is 
a countable transitive model for $ZFC$ (called the {\it ground model}), 
and $\PP \in {\cal M}$. A filter
${\cal G} \subseteq \PP$ is called {\it $\PP$-generic} over ${\cal M}$
iff for all dense $D \in {\cal M}, \; 
{\cal G} \cap D \neq \emptyset$. The countability of ${\cal M}$ implies
that there exist always $\PP$-generic ${\cal G}$; also, if $\PP$ is
{\it non-trivial} in the sense that $\forall p \in \PP \; \exists q ,r 
\in \PP \; (q \leq p, r \leq p \; \land \; q \bot r)$, then a
$\PP$-generic ${\cal G}$ cannot lie in ${\cal M}$, and the
{\it generic extension} ${\cal M}[{\cal G}]$ (the smallest countable
transitive model of $ZFC$ containing ${\cal M}$ and ${\cal G}$)
will be strictly larger than ${\cal M}$. -- The properties of ${\cal M}
[{\cal G}]$ can be described inside ${\cal M}$ using the {\it forcing 
relation} ($\forces$) as follows. For any object in ${\cal M} [{\cal G}]$ there
is a {\it $\PP$-name} in ${\cal M}$; we shall use symbols like
$\breve A$, $\dot A$, ... to denote such names. A sentence of the
{\it forcing language} is a $ZFC$-formula $\psi$ with all free
variables replaced by names. For such $\psi$ and $p \in \PP$ 
we write $p \forces_\PP \;\psi$ ($p$ {\it forces} $\psi$) iff 
for all ${\cal G}$ which are $\PP$-generic over ${\cal M}$, if
$p \in {\cal G}$, then $\psi$ is true in ${\cal M} [{\cal G}]$.
The relation $\forces$ is definable in the ground model ${\cal M}$.
Furthermore, if ${\cal G}$ is $\PP$-generic over ${\cal M}$ and
$\psi$ is true in ${\cal M}[{\cal G}]$, then for some $p \in {\cal G}$,
$p \forces_\PP \;\psi$. \par
An {\it antichain} in a p.o. $\PP$ is a pairwise incompatible set.
$\PP$ is said to satisfy the $\kappa$-cc ({\it $\kappa$-chain condition},
$\kappa$ an uncountable cardinal) iff every antichain $A \subseteq \PP$
has size $< \kappa$. ccc (countable chain condition) is the same as
$\omega_1$-cc. $\PP$ is $\kappa$-closed iff whenever $\lambda < \kappa$
and $\{ p_\xi ; \; \xi < \lambda \}$ is a decreasing sequence of elements
in $\PP$ (i.e. $\xi < \eta \Rightarrow p_\xi \geq p_\eta $), then
$\exists q \in \PP \; \forall \xi < \lambda \; ( q \leq p_\xi)$.
A p.o. $\PP$ {\it preserves cardinals} $\geq \kappa$ ($\leq \kappa$)
iff whenever ${\cal G}$ is $\PP$-generic over ${\cal M}$, and
$\lambda \geq \kappa$ ($\lambda \leq \kappa$, respectively)
is a cardinal in the sense of ${\cal M}$, it is also a cardinal of
${\cal M}[{\cal G}]$. Cardinals which are not preserved are
{\it collapsed}. If $\PP$ has the $\kappa$-cc, then it preserves
cardinals $\geq \kappa$, if it is $\kappa$-closed, it preserves cardinals
$\leq \kappa$. \par
A map $e : \PP \to \QQ$ (where $\PP$ and $\QQ$ are p.o.) is a {\it dense
embedding} iff $\forall p, p' \in \PP \; (p' \leq p \Rightarrow
e(p') \leq e(p))$, $\forall p , p' \in \PP \; (p \bot p' \Rightarrow
e(p) \bot e(p'))$, and $e(\PP)$ is dense in $\QQ$. If $e: \PP \to \QQ$
is dense, $\PP$ and $\QQ$ are equivalent in the sense that they determine
the same generic extensions. Any p.o. can be embedded densely in
a (unique) complete Boolean algebra $\BB (\PP)$ (the Boolean algebra
{\it associated with} $\PP$). \par
Sometimes we want to repeat the generic extension process. This leads to 
the technique of {\it iterated forcing} (see [1] or [8, chapter 2]
for details).
We are mainly concerned with {\it two-step iterations} which we shall
denote by $\PP * \breve \QQ$. \par
We set $Fn(\kappa, \lambda ,\mu) := \{ p ; \; p$ is a function, 
$\vert p \vert < \mu$, $dom(p) \subseteq \kappa$, $ran(p)
\subseteq \lambda \}$; $Fn(\kappa, \lambda, \mu)$ is ordered
by $p \leq q$ iff $p \supseteq q$. $Fn(\kappa, 2, \lambda)$ is called
the ordering for adding $\kappa$ {\it Cohen subsets} of $\lambda$;
for $\lambda = \omega$, the Cohen subsets are referred to as {\it Cohen
reals}. Assume that $2^{<\lambda} = \lambda$, $\lambda$ regular;
then $Fn(\kappa, 2, \lambda)$ is $\lambda$-closed, has the $\lambda^+$-cc,
and so preserves cardinals. Furthermore, if $\kappa^\lambda = \kappa$ (in
the ground model),
then $2^\lambda = \kappa$ in the generic extension.
Cohen extensions can be split and thought of as a two-step iteration
(cf [9, chapter VIII, Theorem 2.1] for the case $\lambda = \omega$). \par
For simplicity, we think of forcing as taking place over the
whole universe ${\cal V}$ instead of over a countable model ${\cal M}$ (though
this is not correct from the formal point of view -- see [9] for
a discussion of this).
\par
Finally we come to {\it internal forcing axioms}. Those are combinatorial
principles proved consistent via iterated forcing; their statement 
captures much of this iteration. The easiest is {\it Martin's Axiom}
$MA$. \smallskip
{\it \item{(MA)} For all ccc p.o. $\PP$ and any family ${\cal D}$ of
$< 2^\omega$ dense subsets of $\PP$, there is a filter ${\cal G}$
in $\PP$ such that $\forall D \in {\cal D} \; ( {\cal G} \cap D
\ne \emptyset)$. \par}
\smallskip
\noindent For the (rather involved) statement of the {\it proper forcing
axiom} $PFA$ we refer the reader to [2] or [8, chapter 3]. 
\bigskip
{\sanse Forcing and inner models.} Sometimes the consistency of $ZFC$ 
is not sufficient for proving the consistency of some combinatorial
statement ($C$) via forcing, and one has to start with a stronger theory
(in general some large cardinal assumption) -- e.g. the existence
of an inaccessible ($ZFC+I$). In those cases we also want to show that
the large cardinal assumption was really necessary; e.g. that
$Con(ZFC+C)$ implies $Con(ZFC+I)$. The way this is usually done is by
showing that if $C$ holds in the universe ${\cal V}$, then some
cardinal is large (e.g. inaccessible) in a {\it sub-universe} ${\cal U}$
(a transitive class model ${\cal U} \subseteq {\cal V}$ satisfying $ZFC$); such
sub-universes are called {\it inner models}. The most important is
the {\it constructible universe} ${\cal L}$, invented by G\"odel.
\par
To show the consistency of the non-existence of weak Kurepa trees,
an inaccessible is collapsed to $\omega_2$ (more correctly, the
cardinals between $\omega_1$ and the inaccessible are collapsed)
-- see [11] or [1]. On the other hand, the non-existence of
Kurepa trees in ${\cal V}$ implies that $\omega_2$ is an inaccessible
cardinal in the sense of ${\cal L}$ (see [9, chapter VII, exercise
(B9)]). The consistency of the existence of Kurepa trees can be
proved by forcing or by showing that they exist in ${\cal L}$.
\vfill\eject

\noindent{\dunhg $\S$ 2. The invariant g(G)}
\Smallskip

{\sanse 2.1.}
For any group $G$ let $g(G)$ -- the
{\it generating number} -- denote
the minimum number of abelian subgroups of $G$ needed to generate
$G$. The following result should be thought of as generalizing the
fact that any countable extraspecial $p$-group is a central sum
of extraspecial $p$-groups of order $p^3$ [16, Corollary 3.10]
-- and so can be generated by two abelian subgroups.
\bigskip
{\sanse 2.2.} {\capit Theorem.} 
{\it For any countable finite-by-abelian group $G$,
$g(G) < \omega$.}
\smallskip
{\it Proof.} We make induction on $\vert G' \vert$. The case $\vert G'
\vert =1$ is trivial. So suppose $\vert G' \vert > 1$. We set
$H:=C_G(G')$. As $G$ is an $FC$-group, $\vert G : H \vert < \omega$;
so it suffices to show that $H$ is generated by finitely many abelian
subgroups. $H'$ is a finite abelian group; i.e. it
is a direct sum of finite cyclic groups of prime power order: 
$H' = \langle a_0 \rangle \oplus ...
\oplus \langle a_n \rangle$. There is a prime $p$ and a natural number
$\ell$ such that $o(a_n)=p^\ell$. Let $A:= \langle a_0 , ... , a_{n-1},
a_n^p \rangle$. We shall define (recursively) two subgroups $H_0 , H_1
\leq  H$ such that $\langle H_0 , H_1 \rangle = H$ and $H'_k \leq A <
H'$ for $k \in 2$. Then the result follows by induction. \par
Suppose $H=\{ b_m ; \; m \in \omega \}$. Let $m_0$ be minimal with the
property that there is an $m$ such that $[ b_{m_0} , b_m ] \not\in A$.
Put $b_0, ..., b_{m_0}$ into $H_0$. Let $c_0 := b_{m_0}$ and $c_1 := b_{m_1}$
where $m_1$ is minimal such that $[c_0, b_{m_1}] \not\in A$, and put $c_1$
into $H_1$. For $m > m_0, m \neq m_1$ let $d_m^0$ be a product of $b_m$
and powers of $c_0$ and $c_1$ such that $[d_m^0 , b_k ] \in A$ for any
$k \in m_0+1 \cup \{ m_1 \}$. We continue this construction recursively.
Suppose we are at step $i$; i.e. $m_{2i}, m_{2i + 1}, c_{2i}, c_{2i + 1}$
and $d_m^i$ ($m > m_{2i}, m \neq m_{2j+1}$ for $j \leq i$) have been 
defined. Then let $m_{2i + 2}$ be minimal with the property that
there is an $m$ such that $[d_{m_{2i + 2}}^i , d_m^i ] \not\in A$.
Put $d_{m_{2i} + 1}^i , ... , d_{m_{2i+2}}^i$ into $H_0$. Let $c_{2i+2}:=
d_{m_{2i+2}}^i$ and $c_{2i + 3} := d_{m_{2i+3}}^i$ where $m_{2i+3}$ is
minimal such that $[c_{2i+2}, d_{m_{2i + 3}}^i ] \not\in A$, and
put $c_{2i + 3}$ into $H_1$. For $m > m_{2i + 2}, m \neq m_{2j+1}$ for
$j \leq i+1$, let $d_m^{i+1}$ be a product of $d_m^i$ and powers of 
$c_{2i + 2}, c_{2i + 3}$ such that $[d_m^{i+1}, d_k^i ] \in A$ for any
$k \in ( (m_{2i}, m_{2i + 2} ] \cup \{ m_{2i + 3} \} ) - \{ m_{2j+1} ;
\; j \leq i \}$.   \par
In the end $H_1 := \langle c_{2j + 1} ; \; j \in \omega \rangle$; and
$H_0$ is the group
generated by the elements which have been put into $H_0$. It is easy to 
see that $H_0$ and $H_1$ satisfy the requirements.  $\qed$
\bigskip
({\it Remark.} The proof of this result is in two steps. The first
shows that finite-by-abelian groups are nilpotent of class 2-by-finite,
and doesn't require countability.)
\bigskip
This property of countable finite-by-abelian groups should
be seen as corresponding to an old result of Baer's, that a group $G$
is centre-by-finite iff $\chi(G) < \omega$ [16, Theorem 7.4],
where $\chi (G)$ denotes the minimum number of abelian subgroups
needed to cover $G$. Nevertheless there are two
drawbacks. First of all it is easy to construct a (countable) $FC$-group
$G$ with $\vert G' \vert = \omega$ but $g(G) = 2$. Secondly, our result
doesn't generalize to higher cardinalities. The important example
of Shelah and Stepr\=ans [13] shows that there are finite-by-abelian
(even extraspecial) groups of size $\omega_1$ all of whose abelian
subgroups are countable. But even for nicer classes of groups there
is nothing corresponding to the Theorem as is shown by the following
\bigskip
{\sanse 2.3.} {\it Example.} Let $E$ be the group 
generated by elements $a, a_\alpha,
\; \alpha < \omega_1$, satisfying the relations $a^p = a_\alpha^p = [a,
a_\alpha] = 1$ and  $[a_\alpha, a_\beta] = a $ for $\alpha < \beta$.
$E$ is easily seen to be an extraspecial ${\cal Z}$-group of exponent $p$.
We will show that $g(E) = \omega_1$.
\par
For suppose that $g(E) \leq \omega$. Then there are abelian subgroups
$A_n$ $(n \in \omega)$ such that $E$ is generated by the $A_n$.
Choose $\Gamma \in 
[\omega_1]^{\omega_1}$ and $n \in \omega$ such that for all $\alpha \in \Gamma$
$a_\alpha \in \langle A_k; \; k < n \rangle$. For each such $\alpha$ and any
$k < n$
we can find $b_{k,\alpha} \in A_k$ such that $a_\alpha = \prod_{k=0}^{n-1}
b_{k,\alpha}$ (at least modulo a factor which is a power of $a$ and
which is irrelevant for our calculation). Now let $B_{k,\alpha}$ consist of the $\beta$ so that 
$a_\beta$
appears as a factor in $b_{k,\alpha}$. 
We may assume that the $B_{k,\alpha}$ form a delta-system with root
$R_k$ for any fixed $k$. Let $C_{k,\alpha} := B_{k,\alpha} - R_k$.
We can suppose that there is a $j_k$ such that $\vert C_{k,\alpha} \vert
= j_k$, that for all $\alpha \in \Gamma$ $sup \; R_k < min \; C_{k,\alpha}$,
that for $\alpha < \beta$ (both in $\Gamma$) 
$sup \; C_{k,\alpha} < min \; C_{k,\beta}$, and that the multiplicities
with which the $a_\beta$ appear in the $b_{k,\alpha}$ depend only on
$\gamma \in R_k$ or $i \in j_k$ (and not on the specific $\alpha$). Then 
$$b_{k,\alpha} = \prod_{\beta \in R_k}  a_\beta^{\ell_\beta} \;
\prod_{i \in j_k}  a_{C_{k,\alpha} (i)}^{m_i},$$
where $\ell_\beta, \; m_i \in p$. An easy commutator calculation
shows that the commutativity of $A_k$ implies
that $\sum_{i \in j_k} m_i \equiv O \; (mod \; p)$. On the other hand,
$$a_\alpha = \prod_{k=0}^{n-1} b_{k,\alpha} = \prod_{k=0}^{n-1}
(\prod_{\beta \in R_k} a_\beta^{\ell_\beta} \;
\prod_{i \in j} a_{C_{k,\alpha} (i)}^{m_i}).$$
This equation cannot hold for any $\alpha$ with $(\{\alpha\}
\cup \cup_{k<n} C_{k,\alpha}) \cap (\cup_{k<n} R_k) = \emptyset$,
thus giving a contradiction. $\qed$                     
\bigskip
{\it Note.} It is easy to see that $E$ can be embedded in an
extraspecial $p$-group
$F$ with $g(F) = 2$. Namely, let $F$ be the group generated by $a$,
$a_\alpha$, $b_\alpha$ ($\alpha < \omega_1$) satisfying -- in addition
to the above relations -- $b_\alpha^p = [a,b_\alpha] = [b_\alpha,b_\beta]
=1$ and
$$[a_\alpha, b_\beta] = \cases{a^{-1} & if $\alpha < \beta$, \cr
1 & otherwise.\cr}$$
Then $F = \langle A_0, A_1 \rangle$, where $A_0= \langle a_\alpha b_\alpha ;
\; \alpha < \omega_1 \rangle$ and $A_1 = \langle b_\alpha ; \;
\alpha < \omega_1 \rangle$. In fact, $F$ is a semidirect extension of $E$. 
\par
So the inequality $g(G) \leq \kappa$ is {\it not} necessarily
preserved when taking subgroups. It is preserved, however, when taking
factor groups. This suggests that instead of dealing with $g$, one
should consider the {\it hereditary generating number} $hg(G) :=
sup \{ g(U) ; \; U \leq G \}$. \par
(A much easier example for this is the direct sum $E_\omega$ of countably 
many extraspecial $p$-groups of size $p^3$ (of exponent $p$ for $p > 2$).
$g(E_\omega) =2$, but $E_\omega$ contains the tree group $C$ of 4.2
which has $g(C) = \omega$.)
\bigskip
{\sanse 2.4.} Let ${\cal QSDF}$ be the ${\cal QSD}$-closure of the class
of finite groups; i.e. $G \in {\cal QSDF}$ iff it is a factor group
of a subgroup of a direct sum of finite groups. ${\cal QSDF}$
is a subclass of ${\cal Z}$ [16, Lemma 3.7]. Tomkinson asked
[16, Question 3F] whether ${\cal Z} \ne {\cal QSDF}$. This was shown
to be true rather indirectly by Tomkinson and L. A. Kurda\v cenko;
namely Kurda\v cenko [10, Theorem 4] proved that any extraspecial
${\cal QSDF}$-group can be embedded in a direct sum of groups of order
$p^3$ with amalgamated centre, and Tomkinson gave a (rather complicated)
example [16, Example 3.16] for an extraspecial ${\cal Z}$-group
which cannot be embedded in a direct sum of groups of order
$p^3$ with amalgamated centre. 
\par
We shall show that the group $E$ of 2.3 does not lie in ${\cal QSDF}$,
thus providing an easier example. To this end, for any group $G$, let 
$P(G)$ be the least cardinal $\kappa$ such that any set of pairwise 
non-commuting elements of $G$ has size less than $\kappa$. A canonical
$\Delta$-system argument shows that $G \in {\cal QSDF}$ implies
$P(G) \leq \omega_1$ (this is a special instance of [3, Theorem 6]).
On the other hand, the definition of $E$ in 2.3 shows that $P(E)=\omega_2$.
Hence $E \in {\cal Z} \setminus {\cal QSDF}$.
\Bigskip
\vfill\eject

\noindent{\dunhg $\S$ 3. FC-automorphisms of countable periodic abelian groups}
\Smallskip

{\sanse 3.1.}
Let $G$ be an $FC$-group. An automorphism $\phi$ of $G$ is called
{\it $FC$-automorphism} iff $\vert \{ g \; (g^{-1}) \phi ; \; g \in G \} \vert
< \omega$; i.e. iff the semidirect extension of $G$ by the group
generated by $\phi$ is still an $FC$-group.
For our discussion the following is important.
\bigskip
{\sanse 3.2.}
{\capit Theorem.} {\it Let $A$ be a countable periodic abelian group.
Suppose $\Phi$ is a group of $FC$-automorphisms of $A$ with $cf(\vert
\Phi \vert) > \omega$. Then there
is a subgroup $B \leq A$ such that $\vert \{B \phi  ; \; \phi \in \Phi \}
\vert = \vert \Phi \vert$.}
\smallskip
{\it Proof} (Tomkinson). First of all, for $\phi \in \Phi$, let $A_\phi :=
\langle a - a \phi ; \; a \in A \rangle$. There are only countably
many finite subgroups $C \leq A$. If $\Phi_C := \{ \phi \in \Phi ;
\; A_\phi \leq C \}$, then $\Phi = \bigcup_{C} \Phi_C$. Since 
$cf(\vert \Phi \vert ) > \omega$, there is a $C=C(A) \leq A$ such
that $\vert \Phi_C \vert = \vert \Phi \vert$.
We make induction on $\vert C \vert$. \par
Secondly we can restrict
our attention to $p$-groups (for some fixed prime $p$). 
The general result follows easily (as any
periodic abelian group is the direct sum of its $p$-components which are
characteristic subgroups). \par
Now let $m := $ exponent of $C$; i.e. $m$ is the smallest integer
such that $p^m C = 0$. Then for all $a \in A$ of height $\geq m$
and all $\phi \in \Phi$, $a \phi = a$. (To see this let $a \in A$ be 
of height $\geq m$. Choose $\hat a \in A$ such that $p^m \hat a = a$.
Let $\hat b := \hat a \phi - \hat a \in C$. Then $a \phi = (p^m
\hat a) \phi = p^m (\hat a + \hat b) = p^m \hat a = a$.)
Especially it suffices to consider reduced $p$-groups. \par
As usual
let $A^1$ denote the subgroup of all elements of infinite height in $A$.
Pr\"ufer's Theorem [5, Theorem 17.3] says that
$A/A^1$ is a direct sum of cyclic $p$-groups.
By the preceding paragraph, $\phi \restrict A^1 = id$ for any
$\phi \in \Phi$. Suppose $A^1 \cap C < C$. Choose $B < A$ containing
$A^1$ such that $B/A^1 \cap (C + A^1)/A^1 = 0$ and $[A:B] < \omega$ (this is
possible because $A/A^1$ is a direct sum of finite groups).
Then either $B$ satisfies the requirements of the Theorem, or $B$ has as
many automorphisms as $A$. In the latter case we are done by induction
because $C(B) < C(A) =C$.
\par
This shows that we may assume $C \leq A^1$ (in particular, $A^1 \neq 0$).
Now let $D < C$ such that $\vert C/D \vert = p$. Each $\phi \in \Phi$
leaves $D$ fixed and so induces an automorphism of $A/D$.
Let $\Phi_{A/D}$ be the group of induced automorphisms. If $\vert
\Phi_{A/D} \vert < \vert \Phi \vert$, then $\vert \Phi_D \vert = \vert
\Phi \vert$,
and we are done by induction. \par
So we may assume that $\vert \Phi_{A/D} \vert = \vert \Phi \vert$ and
consider $\bar A = A/D$. There is an $E \leq A$ such that $E \cap C
= D$ and $A/E \cong C_{p^\infty}$ (such an $E$ can be constructed
as follows: let $\bar F$ be a complement of $\bar C$ in $\{ \bar x
\in \bar A ; \; o(\bar x) = p \}$; set $\bar E := \{ \bar x \in \bar A ;$
if $o(\bar x) = p^n$ then $p^{n-1} \bar x \in \bar F \}$;
let $E$ be the subgroup of $A$ corresponding to $\bar E$). 
For each $\bar\phi \in \Phi_{A/D}$,
$id \restrict p \bar A = \bar \phi \restrict p \bar A$. So $\bar A
= \bar E + p \bar A$ implies that if $\bar \phi  \neq \bar\psi \;
(\bar \phi , \bar\psi \in \Phi_{A/D} )$ then $\bar\phi \restrict
\bar E \neq \bar \psi \restrict \bar E$. Hence $\bar E \cap \bar C 
= \bar 0$ gives us $\bar E \bar \phi \neq \bar E \bar \psi$.
Therefore $E$ has $\vert \Phi_{A/D} \vert$ images under $\Phi$.
This proves the Theorem. $\qed$
\bigskip
{\sanse 3.3.} {\it Proof of Theorems C and C'.}
We have to show that for any finite-by-abelian group $G$ of size
$\omega_1$, \par
\item{(i)} $G \in {\cal Z}$ iff $G \in {\cal Y}$; \par
\item{(ii)} $\vert G / Z(G) \vert \leq \omega$ iff $[G : N_G(A)]
\leq \omega$ for all abelian $A \leq G$. \par
\noindent For suppose not. Then there is a finite-by-abelian group
$G$ which is not ${\cal Z}$ such that \par
\item{} in case (i): $G \in {\cal Y}$; \par
\item{} in case (ii): $[G : N_G (A) ] \leq \omega$ for all abelian $A
\leq G$. \par
\noindent (In case (ii), the fact that $G$ is not in ${\cal Z}$
follows from Tomkinson's result (II) mentioned in the Introduction.)
Then $G$ has a countable
subgroup $U$ with $[G:C_G(U)] = \omega_1$. Let $V := U^G := \langle
g^{-1} U g ; \; g \in G \rangle$. As $G$ is $FC$, $V \nt G$ is
countable. By Theorem 2.2, $g(V) < \omega$, so there
are $n \in \omega$ and $A_i \leq V$ abelian such that $\langle
A_i ; \; i < n \rangle = V$. Clearly $C_G(V) = \cap_{i < n} C_G(A_i)$;
thus there is an $i \in n$ such that $[G:C_G(A_i)] = \omega_1$.
So either $[G:N_G(A_i)] = \omega_1$ in which case we're done, or
$[N_G(A_i) : C_G(A_i)] = \omega_1$. In that case, we may assume
$G=N_G(A_i)$, and $ G / C_G(A_i)$
can be thought of as a group of $FC$-automorphisms of $A_i$, and we are
in the situation of Theorem 3.2; i.e. we get a subgroup $B \leq A_i$
such that $[G : N_G(B)] = \omega_1$, a contradiction. $\qed$
\bigskip
{\sanse 3.4.} The argument in 3.3 shows that if one could 
prove the analogue of Theorem 3.2
under the weaker assumption that $A$ is $FC$ instead of abelian, this
would solve Questions 1 and 1' in the Introduction. So we should ask
\smallskip
{\capit Question 1''.} {\it Suppose $G$ is a countable periodic $FC$-group,
and $\Phi$ is a group of $FC$-automorphisms of $G$ with $cf( \vert \Phi
\vert) > \omega$. Is there a subgroup $U \leq G$ such that
$\vert \{ U \phi ; \; \phi \in \Phi \} \vert = \vert \Phi \vert$ ?}
\vfill\eject

\noindent{\dunhg $\S$ 4. Maximal abelian subgroups of FC-groups}
\Smallskip

{\sanse 4.1.} Maximal abelian subgroups (of $FC$-groups) are important for
our discussion, especially those of countable periodic $FC$-groups 
in case $\kappa = \omega_1$. For Theorem A, we
want to construct a countable $FC$-group having an uncountable
set of automorphisms such that on each (maximal) abelian subgroup only
countably many act differently (4.2 and 4.3). To prove Theorem B, we shall try
to shoot a {\it new} abelian subgroup through an {\it old} set of automorphisms
so that many of these automorphisms act differently on this group (4.4 and 4.5).
These two procedures should be seen as being dual to each other
(cf especially 4.6).
Therefore we pause for an instant to look at the lattice of abelian
subgroups itself.
\bigskip
{\capit Lemma.} {\it If $G$ is a ${\cal Z}_\kappa$-group with
$\vert G/Z(G) \vert \geq \kappa$, then $G$ has at least $2^\kappa$
maximal abelian subgroups.}
\smallskip
{\it Proof.} We construct recursively a tree $\{ A_\sigma ; \; \sigma
\in 2^{<\kappa} \}$ of subgroups of $G$ with $Z(G) \leq A_\sigma$
and $\vert A_\sigma / Z(G) \vert < \kappa$ ($A_\sigma / Z(G)$ is
finitely generated in case $\kappa = \omega$) as follows. Let
$A_{\langle \rangle} := Z(G)$. If $\alpha \in \kappa$ is a limit
ordinal and $\sigma \in 2^\alpha$, let $A_\sigma := \bigcup_{\beta
\in \alpha} A_{\sigma \restrict \beta}$. So assume $\alpha = \beta + 1$
for some $\beta \in \kappa$ and $\sigma \in 2^\beta$. Suppose 
$C_G(A_\sigma)$ is abelian. Choose $B \leq G$ such that $\langle 
B, C_G(A_\sigma) \rangle = G$
and $\vert B \vert < \kappa$ (or $B$ is finitely generated if $\kappa
= \omega$). Then $[G : C_G(B) ] < \kappa$. So $[G : C_G(B) \cap
C_G(A_\sigma)] < \kappa$ which contradicts $\vert G / Z(G) \vert 
\geq \kappa$. So $C_G(A_\sigma)$ is non-abelian and there are $g, \;
h \in C_G(A_\sigma)$ such that $[g,h] \neq 1$. Then set $A_{\sigma
\hat{} \langle 0 \rangle} := \langle A_\sigma, g \rangle$ and
$A_{\sigma \hat{} \langle 1 \rangle} := \langle A_\sigma , h \rangle$.
\par
In the end, for each $f \in 2^\kappa$, extend $\bigcup_{\sigma \subset f}
A_\sigma$ to a maximal abelian subgroup $A_f$. By construction,
$g \neq f$ implies $A_g \neq A_f$. $\qed$ 
\bigskip
In fact, the proof of the Lemma shows that any abelian subgroup $A$
with $\vert AZ(G) / Z(G) \vert < \kappa$ is contained in at least 
$2^\kappa$ distinct maximal abelian subgroups; and that it is contained 
in at least $\kappa$ subgroups $B_\alpha$, $\alpha < \kappa$, with $Z(G)
\leq B_\alpha$ and $\vert B_\alpha / Z(G) \vert < \kappa$ and which are
pairwise {\it incompatible} in the sense that $\langle B_\alpha, B_\beta
\rangle$ is not abelian for $\alpha \neq \beta$ -- this fact will be
used in the proof of Theorem 4.4. below!
\par
As a consequence in case $\kappa = \omega$ we get 
\bigskip
{\capit Corollary.} {\it An $FC$-group has either finitely many or 
at least $2^\omega$ maximal abelian subgroups. It has finitely many
iff it is centre-by-finite.} $\qed$
\bigskip {\sanse 4.2.}
As mentioned earlier we are concerned with the following problem. 
Suppose $G$ is an $FC$-group.
Under which circumstances
is there a set $S$ of $\kappa$ automorphisms of $G$ such that for all
abelian $A \leq G$, $\vert\{ \phi \restrict A ; \; \phi \in S \} \vert <
\kappa$? An easy necessary condition is $g(G) \geq \omega$. We begin
with the following
\bigskip
{\it Example.} For each $n \in \omega$ we introduce a finite group $C_n$
as follows. Let $A_n$ be an elementary abelian $p$-group
of size $p^n$, and $B_n$ an elementary abelian $p$-group of
size $p^{n \choose 2}$. We extend $B_n$ by $A_n$ with factor system
$\tau_n$ as follows:
$$\tau_n(a_i, a_j) = \cases{0 &if $i \geq j$, \cr
b_{h(i,j)} &otherwise,}$$
where $h: [n]^2 \to {n \choose 2}$ is a bijection and the $a_i$ ($b_j$,
respectively) are generators of $A_n$ ($B_n$). Let $C_n$ be the extension
(i.e. $C_n =E(\tau_n)$). Note that $C_n$ is the free object on $n$
generators in the variety of two-step nilpotent groups of exponent
$p$ ($p > 2$); and that it is a
special $p$-group with $C_n' = \Phi(C_n) = Z(C_n) = B_n$.
Let $C$ be the direct sum of the $C_n$. If $g$ is any function from
$\omega$ to $\cup A_n$ with $g(n) \in A_n$, then $g$ defines
a maximal abelian subgroup $M_g := \langle B_n , g(n) ; \; n \in \omega
\rangle$. On the other hand each maximal abelian subgroup of $C$ is
of this form. So the maximal abelian subgroups can be thought of as
branches through a tree. For later reference we shall therefore call $C$
the {\it tree group}. 
\par
Now assume $CH$. Let $\{ M_\alpha ; \; \alpha < \omega_1 \}$ be an 
enumeration of the $M_g$. We introduce (recursively) a set of 
automorphisms $\{ \phi_\alpha ; \; \alpha < \omega_1 \}$ of $G:=C \oplus 
D$ where $D = \langle d \rangle$ 
is a group of order $p$ as follows: fix $\alpha$; let
$\{ N_n ; \; n \in \omega \}$ be an enumeration of $\{ M_\beta ; \;
\beta < \alpha \}$; and let $\{ \psi_n ; \; n \in \omega \}$ be an
enumeration of $\{ \phi_\beta ; \; \beta <  \alpha \}$.
We define $\phi_\alpha$ and an auxiliary function $f: \omega \to \omega$
recursively. Suppose $f \restrict (n+1)$ and 
$\phi_\alpha \restrict (\bigoplus_{i<f(n)} C_i \oplus D)$ have been defined.
We choose $f(n+1)$ so large that we can extend $\phi_\alpha$ to
$(\bigoplus_{i<f(n+1)} C_i \oplus D)$ so that \par
\item{(i)} $\forall c \in \bigoplus_{i < f(n+1)} C_i \oplus D \; \exists
   k \in p \; (c\phi_\alpha  = c + kd)$; \par
\item{(ii)} $\phi_\alpha \restrict (\bigoplus_{i < f(n+1)} 
B_i \oplus D) = id$; \par
\item{(iii)} $\forall k < n \; (\phi_\alpha \restrict 
(\bigoplus_{i < f(n+1)} C_i \oplus D) \neq \psi_k \restrict 
(\bigoplus_{i < f(n+1)} C_i \oplus D))$; \par
\item{(iv)} $\phi_\alpha \restrict ((\cup_{k<n} N_k) \cap 
(\bigoplus_{f(n) \leq i < f(n+1)} C_i \oplus D)) = id$. \par
\noindent This is clearly possible. It is easy to see that $\{ \phi_\alpha
; \; \alpha < \omega_1 \}$ is a set of (distinct) automorphisms of $G$
such that for all maximal abelian $A \leq G$, $\vert 
\{ \phi_\alpha  \restrict A ;
\; \alpha < \omega_1 \} \vert \leq \omega$. $\qed$
\bigskip
{\sanse 4.3.} {\it Proof of Theorem A.}
Form the semidirect extension $E$ of the group $G$ defined in subsection
4.2 by the group $H$ generated
by the automorphisms $\{ \phi_\alpha ; \; \alpha < \omega_1 \}$ (also
defined in 4.2). 
Then $E$ is easily seen to be an $FC$-group with  the required properties
(in fact, for all abelian $A \leq E$, $[E: C_E(A)] \leq \omega$).
$\qed$
\bigskip
{\sanse 4.4.} We now show that $CH$ was necessary in the example above.   
\smallskip
{\capit Theorem.} {\it Let $\lambda > \kappa^+$ be cardinals,
$\kappa$ regular. Denote by $Fn(\lambda,2,\kappa)$ the 
p. o. for adding $\lambda$ Cohen subsets of $\kappa$. Suppose
${\cal V} \models 2^{<\kappa} = \kappa$. Then in ${\cal V} [{\cal G}]$,
where ${\cal G}$ is $Fn(\lambda,2,\kappa)$-generic over ${\cal V}$, the following
holds: for any ${\cal Z}_\kappa$-group $H$ of size $\kappa$ 
and any set of automorphisms
$\Phi$ of $H$ of size $> \kappa$ there is an abelian subgroup $A \leq
H$ such that $\vert\{ \phi \restrict A ; \; \phi \in \Phi \} \vert = \vert
\Phi \vert $.}
\smallskip
{\it Proof.} Let $H$ be any ${\cal Z}_\kappa$-group of size $\kappa$ and
with $\vert H/Z(H) \vert = \kappa$. Note that if $\vert H/Z(H) \vert
< \kappa$, then $Z(H)$ has the required property (as $2^{<\kappa}
=\kappa$).
We define the ordering $\PP_H$ for shooting new abelian subgroups through
$H$ as follows: $\PP_H := \{ A \leq H; \; A $ is abelian and $A$
is generated by
$< \kappa$ elements $\}$. $\PP_H$ is ordered by reverse inclusion; i.e.
$A \leq_{\PP_H} B$ iff $B \subseteq A$. $\PP_H$ is 
$\kappa$-closed, and non-trivial by the discussion in 4.1.
Since $2^{<\kappa} = \kappa$, $\PP_H$ is trivially $\kappa^+ - cc$. 
So forcing with $\PP_H$ preserves cardinals
and cofinalities. 
\par
We first claim that if $\Phi$ is any set of
automorphisms of $H$ of size $> \kappa$ in ${\cal V}$, then in ${\cal V}
[{\cal G}]$, where
${\cal G}$ is $\PP_H$-generic over ${\cal V}$, there is an abelian subgroup $A \leq
H$ such that $\vert\{ \phi \restrict A ; \; \phi \in \Phi \} \vert = \vert
\Phi \vert$.
\par
For $A$ we take the generic object, i.e. $A=\cup\{B \leq H ; \;
B \in {\cal G} \}$. Suppose the claim is false. Let $\mu$ be regular 
with $\vert\Phi\vert \geq \mu \geq \kappa^+$. 
Then there is a $\Psi \subseteq \Phi$ in ${\cal V}[{\cal G}]$ of size
$\mu$ with $\forall \phi, \psi \in \Psi \; (\phi \restrict A
= \psi \restrict A)$. So this statement is forced by a condition
$B \in \PP_H$; i.e. there is a $\PP_H$-name $\breve \Psi$ such that \smallskip
\centerline{$B \forces \breve\Psi \subseteq \Phi \; \land \;
\vert \breve\Psi \vert = \mu \; \land \;
\forall \phi, \psi \in \breve\Psi \;
(\phi \restrict \breve A = \psi \restrict \breve A)$.} \smallskip
\noindent As $\vert \PP_H \vert = 2^{<\kappa} = \kappa$, there is
(in ${\cal V}$) a ${\rm X} \in [\Phi]^{\mu}$ and a $C \leq_{\PP_H} B$
such that \smallskip
\centerline{$C \forces {\rm X} \subseteq \breve\Psi$.}
\smallskip
\noindent Now, $C$ is an abelian
subgroup of the ${\cal Z}_\kappa$-group $H$ of size less than $\kappa$.
So $[H:C_H(C)] < \kappa$. As $\vert {\rm X} \vert > \kappa$ and $2^{<\kappa}
=\kappa$, $\vert \{ \chi \restrict C_H(C) ; \; \chi \in {\rm X}\} \vert
= \mu$ so that we can find $\psi , \chi \in {\rm X}$
and $c \in C_H(C) \setminus C$ such that $\psi (c) \neq \chi (c)$. But then
the condition $\langle C, c \rangle$ forces contradictory statements.
This proves the claim.
\par
Next we remark that for any ${\cal Z}_\kappa$-group $H$ of size $\kappa$,
$\PP_H$ is equivalent (from the forcing theoretic point of view) to
the Cohen forcing $Fn(\kappa,2,\kappa)$ for adding a single new
subset of $\kappa$. For $\kappa = \omega$ this follows from the 
fact that any
two non-trivial countable notions of forcing are equivalent 
[9, chapter VII, exercise (C4), p. 242]. The proof for this 
generalizes as follows. Let $\{ A_\alpha ; \; \alpha < \kappa \}$
enumerate $\PP_H$. We construct recursively a dense embedding $e$
>from $\{ p \in Fn(\kappa, \kappa, \kappa) ; \; dom (p) \in
\kappa \}$ into $\PP_H$. Let $\alpha < \kappa$ and suppose $e \restrict
\{ p \in Fn(\kappa, \kappa, \kappa) ; \; dom (p) \in \alpha \}$
has been defined. If $\alpha$ is limit let, for any $p$ with $dom (p)
= \alpha$, $e(p) = \bigcup_{\beta < \alpha} e(p \restrict \beta)$.
So suppose $\alpha = \beta + 1$ for some $\beta \in \kappa$. There is
by induction
(at least) one $p_0 \in Fn(\kappa, \kappa, \kappa)$ with $dom(p_0) = \beta$
such that $A_\beta$ is compatible with $e(p_0)$. For each $p \in
Fn(\kappa, \kappa, \kappa)$ with $dom (p) = \beta$ choose a maximal 
antichain $M_p$ of size $\kappa$ of conditions below $e(p)$ in
$\PP_H$ (the existence of such an antichain is guaranteed by the
discussion in 4.1), such that $\langle A_\beta, e(p_0)\rangle$ 
is a subgroup of some group in $M_{p_0}$.
Let $e \restrict \{ q \in Fn (\kappa, \kappa, \kappa) ; \; 
dom(q) = \alpha$ and $q \restrict \beta = p \}$ be a bijection onto
$M_p$. It is easy to check that $e$ works. The same argument shows that 
$\{ p \in Fn (\kappa, \kappa, \kappa); \; dom(p) \in \kappa \}$ can be
densely embedded into $Fn (\kappa, 2 , \kappa)$. This gives
equivalence (cf $\S$ 1).
\par
Finally we prove the Theorem. Let $H$ and $\Phi$ be as in the statement
of the Theorem. First suppose $\vert \Phi \vert < \lambda$.
Then $\Phi$ is contained in an initial segment of the extension,
and any subset which is Cohen over this initial segment produces
the required $A$ by the above arguments. So suppose $\vert \Phi \vert
\geq \lambda$. In that case we think of the whole extension as a
two-step extension which first adds $\lambda$ and then $\mu$
Cohen subsets of $\kappa$, where $\lambda \geq \mu \geq \kappa^+$ is
regular with $cf(\vert\Phi\vert) \neq \mu$ and $\vert\Phi\vert > \mu$. 
Then there is a subset $\Psi \in [\Phi]^{\vert \Phi \vert}$ which
is contained in an initial segment of the second extension, and 
our argument applies again. $\qed$     
\bigskip
{\it Remark.} Note that in the Theorem, the assumption $\vert H \vert =
\kappa$ may be replaced by $\vert H/Z(H) \vert = \kappa$. The p.o.
$\PP_H$ in the proof contains in that case the abelian subgroups
$A$ with $Z(H) \leq A$ and $\vert A/Z(H) \vert < \kappa$.
\bigskip
{\sanse 4.5.} {\it Proof of Theorems B and B'.} Let ${\cal V} 
\models ZFC$. We show that
in the model obtained by adding $\omega_2$ Cohen reals to ${\cal V}$, \par
\item{(i)} any ${\cal Y}$-group of size $\omega_1$ is a ${\cal Z}$-group; 
\par
\item{(ii)} there is no $FC$-group
$G$ with $\vert G/Z(G) \vert = \omega_1$ but $[G:N_G(A)] \leq \omega$
for all abelian subgroups $A \leq G$. \par
\noindent For suppose not. Then there is an $FC$-group $G$ of size 
$\omega_1$ which is not in ${\cal Z}$ such that \par
\item{} in case (i): $G \in {\cal Y}$; \par
\item{} in case (ii): $[G:N_G(A)] \leq \omega$ for all abelian $A \leq G$.
\par
\noindent We argue as in the proof of Theorems C and C' (subsection 3.3) using
4.4 instead of 2.2: let $U \leq G$ be countable with 
$[G:C_G(U)] = \omega_1$; let $V:= U^G$.
Apply 4.4 (with $\kappa = \omega$, $\lambda = \omega_2$, $H=V$ and
$\Phi = G/C_G(V)$) to get an abelian $A \leq V$ with $[G : C_G(A)] =
\omega_1$. Now finish as in 3.3 with Theorem 3.2. 
$\qed$
\bigskip
{\sanse 4.6.} The proof of Theorem B shows that its statement
follows from $MA + 2^\omega > \omega_1$. Still this is not the right way to
look at the problem from the combinatorial point of view. 
Namely, when iterating Cohen forcing one merely goes through one {\it particular}
ccc p. o., whereas $MA$ asserts that generic objects exist for {\it all}
ccc p. o. -- not only for those which shoot new abelian subgroups
through an $FC$-group $G$ but also for those which shoot a new automorphism
through $G$ (see below). The consequences of this will
become clear in $\S$ 5 (see the difference between Theorems E and 5.8). \par
Also $MA$ is a weakening of $CH$, and many statements which are provable
in $ZFC + CH$ are still provable in $ZFC + MA$ if we replace $\omega$ by
$<2^\omega$. We shall see now that this is the case for our
problem as well.
\bigskip
{\capit Proposition.} {\it Assume $MA$. Then there is an $FC$-group $G$ with
$\vert G/Z(G) \vert = 2^\omega$ but $[G :N_G(A)] < 2^\omega$ for
all abelian subgroups $A \leq G$.}
\smallskip
{\it Proof.} Let $C$ be again the tree group of 4.2. We define
the partial order $\QQ_C$ for shooting new automorphisms through $C
\oplus D$ (where $D = \langle d \rangle$ is again a group of order $p$).
$\QQ_C := \{ (\phi, {\cal A}); \; \phi$ is a finite partial automorphism
of $C \oplus D$ with (i) $\exists n \in \omega$ with 
$dom(\phi) = \bigoplus_{i \in n} C_i \oplus D$, (ii) $\forall c \in 
dom(\phi) \; \exists k \in p
\; ( c\phi = c + kd)$ and (iii) $\phi \restrict  
(\bigoplus_{i \in n} B_i \oplus D ) = id$; and ${\cal A}$
is a finite collection of maximal abelian subgroups of $C \; \}$;
$(\phi, {\cal A}) \leq_{\QQ_C} (\psi, {\cal B})$ iff $\phi \supseteq
\psi$ and ${\cal A} \supseteq {\cal B}$ and $\forall c \in (dom(\phi)
-dom(\psi)) \cap (\cup{\cal B}) \; (c\phi = c)$. $\QQ_C$ is ccc
and generically shoots a {\it new} automorphism through $C \oplus D$ which
equals the identity on all {\it old} abelian subgroups from some point on.
\par
To prove the Proposition let $\{ A_\alpha ; \; \alpha < 2^\omega \}$ enumerate
the maximal abelian subgroups of $C$. We construct recursively
a set of automorphisms $\{ \phi_\alpha ; \; \alpha < 2^\omega \}$.
Let ${\cal D}_\alpha :=
\{ D_{\beta} ; \; \beta  < \alpha \}$ where $D_{\beta} :=
\{ (\phi, {\cal A}) \in \QQ_C ; \; A_\beta \in {\cal A}$ and $c \phi \neq
c \phi_\beta$ for some $c \in dom(\phi) \}$. Each 
$D_{\beta}$ is dense in $\QQ_C$; hence, by $MA$, there is a
${\cal D}_\alpha$-generic filter $G_\alpha$.
Let $\phi_\alpha := \cup \{ \phi; \; \exists {\cal A} \; (\phi,
{\cal A}) \in G_\alpha \}$. Then for all maximal abelian $A \leq C 
\oplus D$, $\vert \{ \phi_\alpha \restrict A ; \; \alpha < 2^\omega \}
\vert < 2^\omega$. Now let $G:=(C \oplus D) \sd \langle \phi_\alpha ; \;
\alpha \in 2^\omega \rangle$. $\qed$
\bigskip
Certainly one should ask whether $MA$ is necessary at all in the above
result; or whether it can be proved in $ZFC$ alone. It turns out
that the answer (to the second question) is no, at least if we assume
the existence of an inaccessible cardinal -- see $\S$ 5 (Theorem E).
\bigskip
{\sanse 4.7.} It is quite usual that combinatorial statements are not
decided by $\lnot CH$. Again this is true in our situation.
\smallskip
{\capit Proposition.} {\it It is consistent that $2^\omega > \omega_1$
and there is an $FC$-group $G$ with $\vert G /Z(G) \vert = \omega_1$,
but $[G : N_G(A) ] \leq \omega$ for all abelian subgroups
$A \leq G$.}
\smallskip
{\it Sketch of the proof.} 
The proof uses the tree group of 4.2 as main ingredient.
Start with ${\cal V} \models 2^\omega > \omega_1$ and make a finite support 
iteration of length $\omega_1$ of the partial order $\QQ_C$
described in 4.6.
$\qed$
\vfill\eject

\noindent{\dunhg $\S$ 5. Extraspecial p-groups and Kurepa trees}
\Smallskip

{\sanse 5.1.} The goal of this section is a detailed investigation of
extraspecial $p$-groups, especially of those of size $\omega_2$.
The philosophy behind this is that many {\it bad} things that
can happen to (periodic) $FC$-groups already happen in case
of extraspecial $p$-groups; or even that {\it bad} periodic
$FC$-groups involve {\it bad} extraspecial groups -- 
the most surprising example for this is Tomkinson's result that a 
periodic $FC$-group $G$ which does not lie in ${\cal Y}$
contains $U \nt V \leq G$ such that $V/U$ is extraspecial and
not in ${\cal Y}$ (see [14] or [16, Theorem 3.15]).
Our main contribution
in this direction is the equiconsistency result mentioned in
the Introduction (Corollary to Theorems D and E). Another example
is the equivalence in Proposition 5.5. -- On the
other hand, because of their simple {\it algebraic} structure (e.g.,
the fact that subgroups are either normal or abelian),
extraspecial examples are in general the easiest to construct,
and such constructions depend only on the underlying {\it
combinatorial} structure -- the classical example for this is the
existence of a Shelah-Stepr\=ans group [13].
\par
We let ${\cal Y}_\kappa$ be the class of groups in which $[G:N_G(U)]
< \kappa$ whenever $U \leq G$ is generated by fewer than $\kappa$
elements. So the class ${\cal Y}$ is just the class of locally finite
groups in the intersection of the ${\cal Y}_\kappa$; and also
${\cal Y}_\omega = {\cal Z}_\omega =$ the class of $FC$-groups.
It follows from Theorems B' and C' that ${\cal Y}_{\omega_1}$ and
${\cal Z}_{\omega_1}$ are consistently equal for periodic
$FC$-groups, and that they are equal for periodic finite-by-abelian
groups.
\par
By Tomkinson's result mentioned in the Introduction (II), if 
$\lambda < \kappa$ are cardinals and $G$ is a group with
$\vert G/Z(G) \vert = \kappa$ and $[G:N_G(U)] \leq \lambda$
for all subgroups $U \leq G$, then $G$ is ${\cal Y}_\mu$ for any
$\mu \geq \lambda^+$ but not ${\cal Z}_\kappa$. Furthermore,
for extraspecial $p$-groups $G$, the following are equivalent
(where $\kappa$ is any cardinal). \par
\item{(i)} $[G:N_G(A)] < \kappa$ for all (abelian) subgroups $A \leq
G$. \par
\item{(ii)} For all maximal abelian subgroups $A$ of $G$, $[G:A]
< \kappa$. \par
\noindent  In particular, an extraspecial $p$-group of size $\omega_2$
whose maximal abelian subgroups satisfy $[G:A] \leq \omega_1$
is ${\cal Y}_{\omega_2}$ but not ${\cal Z}_{\omega_2}$. 
\par
This should motivate us to study the three classes ${\cal Y}_{\omega_1}
= {\cal Z}_{\omega_1}$, ${\cal Y}_{\omega_2}$, and ${\cal Z}_{\omega_2}$
for extraspecial $p$-groups more thoroughly. 
Clearly, there are groups lying in none
or in all of these classes, or in ${\cal Z}_{\omega_2} \setminus
{\cal Z}_{\omega_1}$. The existence of groups which are in
${\cal Y}_{\omega_2} \setminus {\cal Z}_{\omega_2}$ or in
${\cal Z}_{\omega_1} \setminus {\cal Z}_{\omega_2}$ will be discussed
in the subsequent subsections (up to 5.5). Our results can
be summarized in the following chart.
\bigskip
{\offinterlineskip \tabskip=0pt
\halign{ \strut \vrule#& \quad # & \vrule#& \quad # & \vrule#&
\quad # & \vrule# \cr
\noalign{\hrule}
&&&&&& \cr
&   && ${\cal Z}_{\omega_1}={\cal Y}_{\omega_1}$ &&
$\neg {\cal Z}_{\omega_1} = \neg {\cal Y}_{\omega_1}$ & \cr
&&&&&& \cr
\noalign{\hrule}
&&&&&& \cr
& ${\cal Z}_{\omega_2}$ && easy && easy & \cr
&&&&&& \cr
\noalign{\hrule}
&&&&&& \cr
&  &&
 && 5.3. and 5.4. (follows from the & \cr
& $\neg {\cal Z}_{\omega_2}$ but ${\cal Y}_{\omega_2}$
  && ? (cf. 5.5.)  && existence of Kurepa trees, and implies & \cr
&   &&   && the existence of weak Kurepa trees) & \cr
&&&&&& \cr
\noalign{\hrule}
&&&&&& \cr
& $\neg {\cal Y}_{\omega_2}$ && 5.5. (equivalent to the && easy & \cr
& && existence of Kurepa trees) && & \cr
&&&&&& \cr
\noalign{\hrule} }}
\bigskip\bigskip
{\sanse 5.2.} The following is useful for the proof of Theorem D.
\smallskip
{\capit Lemma {\rm (Folklore)}.} 
{\it Assume there is a Kurepa family. Then there
is an a. d. Kurepa family of the same size.}
\smallskip
{\it Proof.} Let $\{ A_\alpha ; \; \alpha < \kappa \}$ be a Kurepa family
(where $\kappa \geq \omega_2$).
Let $f$ be a bijection between $\{ A_\alpha \cap \beta ; \;
\alpha < \kappa , \; \beta < \omega_1\}$ and $\omega_1$. 
Then $\{ \{f(A_\alpha \cap
\beta) ; \; \beta < \omega_1 \} ; \; \alpha < \kappa \}$ is easily
seen to be an a. d. Kurepa family. $\qed$
\bigskip
{\sanse 5.3.} {\it Proof of Theorem D.}
Let $E$ be a Shelah-Stepr\=ans-group [13] of size
$\omega_1$, and let ${\cal A} =\{ A_\alpha ; \; \alpha < \omega_2 \}$ be
an a. d. Kurepa family. We extend $E$ semidirectly by an elementary abelian
group $B$ of automorphisms using ${\cal A}$ as follows:
for all $\alpha < \omega_2$ define $\phi_\alpha$ by
$$a_\beta\phi_\alpha  = \cases{ a_\beta &if $\beta = 0$ or $\beta \not\in
A_\alpha$, \cr
a_\beta  a_0 &otherwise,\cr}$$
\noindent where $a_0$ generates $Z(E)$ and $\{ a_\beta ; \; 1 \leq \beta
< \omega_1 \}$ generates $E$. Set $B:= \langle \phi_\alpha ; \; \alpha
< \omega_2 \rangle$. This completes the construction. $G := E \sd B$
is easily seen to be extraspecial.
\par
Now suppose $A \leq G$ is abelian. Let $\pi(A)$ denote the subgroup
of $E$ generated by 
the projection of $A$ on the first coordinate (we think of the
semidirect product as a set of tuples). We claim that $\pi(A)$ is
countable. For suppose not. Then clearly $a_0 \in \pi(A)$.
Let $C$ be a maximal abelian subgroup of $\pi(A)$. $C$ is countable,
and $C_{\pi(A)}(C)=C$. Choose a subset $\{ (b_\alpha, \psi_\alpha) ;
\; \alpha < \omega_1 \}$ of $A$ such that $C \leq \langle b_n ; \; n
\in \omega \rangle$ and $b_\alpha \neq b_\beta$ for $\alpha \neq
\beta$. Now let $B_{\alpha}$ consist of the $\beta$ so that 
$a_\beta$
appears as a factor in $b_{\alpha}$. 
We may assume that the $B_{\alpha}$ ($\alpha \geq
\omega$) form a delta-system with root
$R$. Let $C_{\alpha} := B_{\alpha} \setminus R$.
We can suppose that there is a $j$ such that $\vert C_{\alpha} \vert
= j$ for $\alpha \geq \omega$, that for $\alpha < \beta$ we have
$sup \; C_{\alpha} < min \; C_{\beta}$, and that the multiplicities
with which the $a_\beta$ appear in the $b_\alpha$ depend only on
$\gamma \in R$ or $i \in j$. 
As ${\cal A}$ is a. d., we may assume that for each of the (countably
many) automorphisms $\phi_\delta$ appearing as a factor in some $\psi_n$
($n \in \omega$) and each $i \in j$ either $\forall \alpha \geq \omega
\; (a_{C_\alpha (i)} \phi_\delta = a_{C_\alpha (i)} )$ or $\forall
\alpha \geq \omega \; (a_{C_\alpha (i)} \phi_\delta = a_{C_\alpha (i)}
a_0 )$ (without loss the corresponding $A_\delta$'s are disjoint above
$C_\omega (0)$). In particular we have that 
for fixed $n \in \omega$,
$c_n := b_\alpha^{-1}(b_\alpha)\psi_n = 
b_\beta^{-1} \; (b_\beta)\psi_n$ for any $\alpha, 
\beta \geq \omega$. As ${\cal A}$ is a Kurepa family, we may assume
that $\psi_\alpha \restrict C = \psi_\beta \restrict C$ for
any $\alpha, \beta  \geq \omega$. But then
$$[(b_n, \psi_n), (b_\alpha , \psi_\alpha)] =
((b_n^{-1}\psi_n^{-1}) \psi_\alpha^{-1} \; (b_\alpha^{-1}) \psi_\alpha^{-1},
\psi_n^{-1} \psi_\alpha^{-1} ) \; (b_n \psi_\alpha \; b_\alpha , \psi_n
\psi_\alpha )$$ $$ = (b_n^{-1} \; (b_\alpha^{-1})\psi_n \; (b_n) \psi_\alpha
\; b_\alpha , 1) = (c_n^{-1} d_n [b_n,b_\alpha] ,1),$$
\noindent where $d_n = b_n \psi_\alpha \; b_n^{-1}$. As $C$ is maximal 
abelian in $\pi(A)$,
there is certainly an $n \in \omega$ such that $[b_n,b_\alpha]
\neq [b_n,b_\beta]$ for some $\alpha, \beta \geq \omega$. But then
the above calculation shows that $A$ cannot be abelian. \par
Now the fact that $\pi(A)$ is countable and that ${\cal A}$ is
a Kurepa family implies $[G:C_G(\pi(A))] \leq \omega_1$ (in fact,
equality holds unless $\pi(A)$ is finite, 
because $E$ is a Shelah-Stepr\=ans-group).
If $\rho(A)$ is the projection of $A$ on the second coordinate,
$[G:C_G(\rho(A))] \leq \omega_1$ holds trivially; and $A \leq \langle
\pi(A), \rho(A) \rangle$ implies $[G:N_G(A)] \leq [G:C_G(A)]
\leq \omega_1$. $\qed$
\bigskip

{\sanse 5.4.} {\capit Theorem.} {\it If there is an extraspecial $p$-group of size
$\omega_2$ in ${\cal Y}_{\omega_2}$ but not in ${\cal Z}_{\omega_2}
$, then there is a weak Kurepa tree.}
\smallskip
{\it Proof.} Let $G$ be such a group. Choose $U \leq G$ of size $\omega_1$
such that $[G:C_G(U)] = \omega_2$. Let $\{ u_\alpha ; \; \alpha <
\omega_1 \}$ generate $U$. Let $\{ f_\beta ; \; \beta < \omega_2 \}$
be a subset of $G \setminus U$ such that $f_\alpha C_G(U) \neq f_\beta C_G(U)$.
Define $g_\beta : \omega_1 \to p$ for $\beta < \omega_2$
by $g_\beta (\alpha) = k$ iff
$[u_\alpha, f_\beta] = ka$, where $a$ generates $G'$.
We claim that the $g_\beta$ form the branches of a weak Kurepa tree.
\par
For suppose not. Then there is an $\alpha \in \omega_1$ such that $\vert
\{ g_\beta \restrict
\alpha ; \; \beta < \omega_2 \} \vert = \omega_2$. This immediately
implies that $[G:C_G(V)] = \omega_2$ for a countable subgroup $V \leq U$.
$V$ is a direct sum of an extraspecial and an abelian group;
especially $g(V) \leq 2$ (this follows from the fact that countable
extraspecial $p$-groups are central sums of groups of order $p^3$).
So there is an abelian $A \leq V$ such that $[G:C_G(A)] = \omega_2$.
Cutting away $G'$ if necessary we may assume that $[G:N_G(A)] =
\omega_2$, contradicting the fact that $G \in {\cal Y}_{\omega_2}$.
$\qed$
\bigskip
Note that in the hypothesis of the Theorem, {\it extraspecial $p$-group}
can be replaced by {\it (periodic) finite-by-abelian group}. To see that this
more general result is true, just apply Theorems 2.2 and 3.2
at the end of the proof. And if Question 1'' had a positive answer,
we could prove this for {\it (periodic) $FC$-groups}. 
\par
There is a gap between Theorem D and Theorem 5.4. We feel
that it should be possible to make a construction like the one in
5.3 using a weak Kurepa tree only.
\bigskip
{\sanse 5.5.} Using the same techniques as in 5.3 and 5.4 we get
\smallskip
{\capit Proposition.} {\it The following are equivalent. \par
\item{(i)} There is a Kurepa tree. \par
\item{(ii)} There is an extraspecial $p$-group which is 
${\cal Z}_{\omega_1}$ but not ${\cal Z}_{\omega_2}$. \par
\item{(iii)} There is an $FC$-group which is 
${\cal Z}_{\omega_1}$ but not ${\cal Z}_{\omega_2}$.}
\smallskip
{\it Proof.} To see one direction ( $(i) \Rightarrow (ii)$ ) let $E$
be any extraspecial ${\cal Z}$-group of size $\omega_1$, and let, as
in 5.3, ${\cal A} = \{ A_\alpha ; \; \alpha < \omega_2 \}$ be an
a.d. Kurepa family. For all $\alpha < \omega_2$
define $\phi_\alpha$ by
$$a_\beta \phi_\alpha = \cases{ a_\beta &if $\beta = 0$ or $\beta
\not\in A_\alpha$, \cr 
a_\beta a_0 &otherwise, \cr}$$
where $a_0$ generates $Z(E)$ and $\{ a_\beta ; \; 1 \leq \beta < \omega_1
\}$ generates $E$. Set $G:= E \sd \langle \phi_\alpha ; \;
\alpha < \omega_2 \rangle$. Clearly $G$ has the required properties.
\par
Conversely, to see $(iii) \Rightarrow (i)$, make the same construction
as in 5.4. $\qed$
\bigskip
In general, the group constructed in the first part of the proof
will not lie in ${\cal Y}_{\omega_2}$ either. Hence the only question
left open is whether there are extraspecial ${\cal Z}_{\omega_1}$-groups
in ${\cal Y}_{\omega_2} \setminus {\cal Z}_{\omega_2}$. We conjecture
that they exist in the constructible universe ${\cal L}$. Such a group
of size $\omega_2$ would lie in ${\cal Y} \setminus {\cal Z}$
as well and so give an answer to Question 3F in [16]. \par
On the other hand, unlike the other classes considered so far,
the consistency of $ZFC$ {\it alone} implies the consistency of the
non-existence of extraspecial ${\cal Z}_{\omega_1}$-groups in
${\cal Y}_{\omega_2} \setminus {\cal Z}_{\omega_2}$. To see this,
let ${\cal V} \models ZFC + GCH$. Add $\omega_3$ Cohen subsets of
$\omega_1$. We claim that in the resulting model ${\cal V} [{\cal G}]$,
there are no such groups. For suppose $G$ is such a group.
Find $U \leq G$ of size $\omega_1$ with $[G:C_G(U)] = \omega_2$,
without loss $U \nt G$. Apply 4.4 with $\lambda = \omega_3$,
$\kappa = \omega_1$, $H=U$, and $\Phi = G/C_G(U)$ (this can be done as
$U \in {\cal Z}_{\omega_1}$). Find $A \leq U$ abelian such that
$[G :C_G(A)] = \omega_2$. Cutting $G'$ away, if necessary,
we can assume $C_G(A) = N_G(A)$, a contradiction.
\bigskip
{\sanse 5.6.} We now want to turn to the proof of Theorem E.
Certainly, in a model where its statement is true, neither of the
{\it bad} situations discussed in 4.2 (and 4.3, 4.6, 4.7) and
in 5.3 can occur. So we'd better look for a model where there
are no Kurepa trees and where $CH$ is false. The discussion in 5.4
and 4.4 suggests that there shouldn't be weak Kurepa trees either
and that there should be reals Cohen over ${\cal L}$. One possible attack
would be to destroy all weak Kurepa trees by collapsing an
inaccessible to $\omega_2$ (as in [11, $\S\S$3,4] or [1, $\S$8])
and then to add $\omega_2$ Cohen reals (by the last $\omega_2$
we mean, of course, the $\omega_2$ of the intermediate model).
Unfortunately, we don't know whether it is true in general that
there are no weak Kurepa trees in the final extension; but this
is true if the intermediate model is Mitchell's [11].
In fact, it turns out that in this case the second extension is 
unnecessary, and what we want to show consistent already holds
in Mitchell's model. The reason for this is essentially that this
model is gotten by first adding $\kappa$ Cohen reals (where $\kappa$
is inaccessible) and then collapsing $\kappa$ to $\omega_2$
using a forcing which does not add reals -- and hence does not
destroy the {\it nice} situation created by the Cohen reals --
while killing all weak Kurepa trees.
\par
First we will review Mitchell's model and some elementary facts
about it. Let ${\cal V} \models $"$ ZFC+GCH+$ there is an inaccessible".
Let $\kappa$ be inaccessible in ${\cal V}$. $\PP = Fn(\kappa, 2, \omega)$
is the ordering for adding $\kappa$ Cohen reals. $\BB = \BB (\PP)$
is the Boolean algebra associated with $\PP$. Set $\PP_\alpha :=
\{ p \in \PP ; \; supp(p) \subseteq \alpha \}$ for $\alpha < \kappa$.
$\BB_\alpha$ is the Boolean algebra associated with $\PP_\alpha$.
$f \in {\cal V}$ is in the set ${\cal A}$ of {\it acceptable functions}
iff \par
\item{(1)} $dom(f) \subseteq \kappa$; $ran(f) \subseteq \BB$; \par
\item{(2)} $\vert dom(f) \vert \leq \omega$; \par
\item{(3)} $f(\gamma) \in \BB_{\gamma+\omega}$ for $\gamma \in dom(f)$.
\par
\noindent If ${\cal F}$ is $\PP$-generic over ${\cal V}$,
$f \in {\cal A}$, then we define $\bar f : dom(f) \to 2$ in ${\cal V}[{\cal F}]$
by $\bar f(\gamma) = 1$ iff $ f(\gamma) \in {\cal F} $.
Define $\QQ$ in ${\cal V}[{\cal F}]$ by letting the underlying
set of $\QQ$ be ${\cal A}$ and $f \leq_{\QQ} g$ iff $\bar f \supseteq \bar g$.
So we get a 2-step iteration $\PP * \QQ$ with $(p,f) \leq (q,g)$ iff
$p \leq q$ and $p \forces_\PP f \leq_{\QQ} g$. We shall denote the final
extension by ${\cal V} [{\cal F}] [{\cal G}]$.
\smallskip
{\capit Facts {\rm (Mitchell [11])}.} {\it (1) Suppose $p \forces_\PP
$"$\breve D$ is open dense in $\QQ$ below $f$", where $f \in {\cal A}$.
Let $g \in {\cal A}$ such that $g \supseteq f$. Then there is $h \in {\cal A}$
such that $h \supseteq g$ and $p \forces_\PP h \in \breve D$.
\par
(2) $\QQ$ does not add new functions with countable domain over
${\cal V}[{\cal F}]$; i.e. if $t : \omega \to {\cal V} [{\cal F}]$
where $t \in {\cal V} [{\cal F}] [{\cal G}]$, then $t \in {\cal V}
[{\cal F}]$. 
\par
(3) Let $\{ g_\alpha ; \; \alpha < \kappa \} \subseteq {\cal A}$. Then there
are $X \in [\kappa]^\kappa$ and $g \in {\cal A}$ such that $\forall \alpha
, \beta \in X \; (dom(g_\alpha) \cap dom(g_\beta) = dom (g))$ 
and $\forall \alpha \in X \; (g_\alpha \restrict dom (g) =g)$. 
\par
(4) $\PP * \QQ$ preserves $\omega_1$ (this follows from the ccc-ness
of $\PP$ and fact (2)) and cardinals $\geq \kappa$ (it follows from (3)
that $\PP * \QQ$ is $\kappa$-cc), but collapses all cardinals
in between to $\omega_1$; i.e. $\kappa^{\cal V} = \omega_2^{{\cal V}
[{\cal F}] [{\cal G}]}$. 
\par
(5) In ${\cal V}[{\cal F}] [{\cal G}]$, $2^\omega = 2^{\omega_1} =
\omega_2$.
\par
(6) Let $\nu < \kappa$ be such that $\nu' + \omega \leq \nu$
for each $\nu' < \nu$. Then the generic extension via $\PP * \QQ$
can be split in a 2-step extension, the first of which adds
$\vert\nu\vert$ Cohen reals and a $\QQ \restrict \nu$-generic function from
$\nu$ to $2$, whereas the second adds the remaining $\kappa$
Cohen reals and the remaining part of the $\QQ$-generic function
>from $\kappa$ to $2$.
\par
(7) In ${\cal V}[{\cal F}] [{\cal G}]$, there are no weak Kurepa trees.}
\smallskip
{\it Proofs.} (1) to (5) are (more or less) 3.1 to 3.5 in [11];
concerning (3) we note that it is proved via a straightforward
$\Delta$-system-argument. (6) is made more explicit on pp.
29 and 30 in [11] and proved in 3.6. For (7), see 4.7
in [11]. $\qed$
\bigskip
{\sanse 5.7.} {\it Proof of Theorem E.} We show that
${\cal V} [{\cal F}] [{\cal G}] \models$"for both $\omega_1$
and $\omega_2$ and any $FC$-group $G$, (i) through (iii) in (II)
are equivalent", where ${\cal V}[{\cal F}][{\cal G}]$
is Mitchell's model as in the preceding section.
\par
Counterexamples $G$ with $\vert G/Z(G) \vert = \omega_1$ are easily
excluded. Without loss such $G$ would have size $\omega_1$. By the 
$\kappa$-cc of $\PP * \QQ$ (which follows from the ccc of $\PP$
and fact (3) in 5.6) it would lie in an intermediate extension
(see fact (6)). Any real Cohen over this intermediate extension
shows that the assumption was false (by the argument of Theorem
4.4).
\par
Suppose $G$ is a counterexample with $\vert G/Z(G) \vert = \omega_2$;
without loss $\vert G \vert = \omega_2$; $[G:N_G(A)] \leq \omega_1$
for all abelian $A \leq G$; and there is a $U \leq G$ of size $
\leq\omega_1$ such that $[G : C_G(U)] = \omega_2$ (because $G$ cannot be
a ${\cal Z}_{\omega_2}$-group by Tomkinson's result (II) in the 
Introduction). Without loss $\vert U \vert = \omega_1$.
Let $V := U^G = \langle x^{-1} U x ; \; x \in G \rangle$.
As $G$ is an $FC$-group, $\vert V \vert = \omega_1$; and $V \nt G$.
Clearly $\vert G/C_G(V) \vert = \omega_2$. For any $\bar g \in
G/C_G(V)$ define a function $f_{\bar g}: V \to V$ by $f_{\bar g}
(v) := g^{-1} v g v^{-1}$ where $g \in \bar g$ is arbitrary.
Think of $\{ f_{\bar g} ; \; \bar g \in G/C_G(V) \}$ as the
set of branches through a tree $T$. As ${\cal V}[{\cal F}][{\cal G}]$
does not contain weak Kurepa trees (fact (7) in 5.6), there is a countable
$S \subseteq V$ such that $\{ f_{\bar g} \restrict S ; \;
\bar G \in G/C_G(V) \}$ has size $\omega_2$. Let $W :=
\langle S^G \rangle$. $W$ is a countable normal subgroup
of $G$; and $\vert G / C_G(W) \vert = \omega_2$ by construction.
\par
We now want to prove that there is an abelian $A \leq W$ such
that $[G:C_G(A)] = \omega_2$ ({\it main claim}). The way we do this is an elaboration of 
the proof of Theorem 4.4. For this argument it is crucial that we use
Mitchell's model and not just any model without weak Kurepa trees.
\par
We think of $G/C_G(W)$ as a group of automorphisms $\Phi$ of $W$;
more explicitly, $\Phi : \omega_2 \to Aut(W)$. In ${\cal V} [{\cal F}]$,
let $\breve \Phi$ be a $\QQ$-name for $\Phi$. Let $D_\alpha$
($\alpha < \kappa$) be the set of conditions deciding $\breve \Phi
(\alpha)$. $D_\alpha$ is open dense by fact (2). Let $\dot D_\alpha$
be a $\PP$-name for $D_\alpha$ ($\alpha < \kappa$). Then
$$\forces_\PP "\dot D_\alpha \; {\rm  is \; open \; dense }".$$
By fact (1) there are $g_\alpha \in {\cal A}$ such that
$$\forces_\PP \;  g_\alpha \in \dot D_\alpha.$$
So in ${\cal V}[{\cal F}]$, there are $\phi_\alpha$ such that
$g_\alpha \forces_{\QQ} \breve \Phi (\alpha) = \phi_\alpha$. Let
$\dot \phi_\alpha$ ($\alpha < \kappa$) be $\PP$-names for the
$\phi_\alpha$. Then $\forces_\PP \; g_\alpha \forces_{\QQ} \dot
\Phi(\alpha) = \dot \phi_\alpha$, where $\dot\Phi$ is a $\PP$-name
for $\breve\Phi$. Using fact (3) we get $X \in [\kappa]^\kappa$
and $g \in {\cal A}$ such that $\forall \alpha, \beta \in X$,
$dom(g_\alpha) \cap dom(g_\beta) = dom (g)$ and $\forall \alpha
\in X$, $g_\alpha \restrict dom(g) = g$. Now we split $\PP$ into two
parts (i.e. $\PP = \PP_1 \times \PP_2$) such that \par
\item{(1)} $\PP_1$ adds $\kappa$ Cohen reals and $\PP_2$ adds
one Cohen real; \par
\item{(2)} there is $Y \in [X]^\kappa$ such that $\forall \alpha \in
Y \; (\phi_\alpha \in {\cal V}[{\cal F}_1])$ where ${\cal F}_1$
is $\PP_1$-generic over ${\cal V}$. \par So in ${\cal V}[{\cal F}_1]$,
$$\forces_{\PP_2} \; g_\alpha \forces_{\QQ} \dot \Phi (\alpha) = \phi_\alpha,$$
where $\alpha \in Y$. From now on we work in ${\cal V}[{\cal F}_1]$.
As in the proof of Theorem 4.4 we think of $\PP_2$ as adding a new abelian
subgroup of $W$. Let $\dot A$ be a $\PP_2$-name for this generic object.
\par
The rest of the proof of the main claim is by contradiction.
Suppose that
$$\forces_{\PP * \QQ} \forall B \leq W \; {\rm  abelian } \;\;
( \vert \{ \dot\Phi (\alpha) \restrict B ; \; \alpha < \kappa \} \vert
< \kappa).$$
Especially, in ${\cal V}[{\cal F}_1]$,
$$\forces_{\PP_2 * \QQ} \vert \{ \dot\Phi (\alpha) \restrict \dot A
; \; \alpha \in Y \} \vert < \kappa.$$
Hence,
$$\forces_{\PP_2 * \QQ} \exists \alpha \; \forall \beta \geq \alpha \;
(\beta \in Y \Rightarrow \exists \gamma < \alpha \;
(\gamma \in Y \wedge \dot\Phi(\beta)\restrict \dot A = \dot\Phi (\gamma)
\restrict \dot A )).$$
So there are $\alpha < \kappa$ and $C \in \PP_2$, $h \supseteq g$,
$h \in {\cal A}$ such that
$$(C,h) \forces_{\PP_2 * \QQ}  \forall \beta \geq \alpha \;
(\beta \in Y \Rightarrow \exists \gamma < \alpha \;
(\gamma \in Y \wedge \dot\Phi(\beta)\restrict \dot A = \dot\Phi (\gamma)
\restrict \dot A )).$$
Choose $Z \in [Y\setminus \alpha]^\kappa$ such that for
$\beta \in Z$, $g_\beta$ and $h$ are compatible (in ${\cal V}$).
For $\beta \in Z$ let (in ${\cal V}[{\cal F}]$) $D_\beta$ be the set
of conditions forcing $\breve \Phi(\beta) \restrict A = \breve
\Phi(\gamma) \restrict A$ for some $\gamma < \alpha$. $D_\beta$
is open dense below $h$. Let $\dot D_\beta$ be a $\PP_2$-name for
$D_\beta$. By fact (1), there are $h_\beta$ such that $h_\beta 
\supseteq g_\beta$ and $h_\beta \supseteq h$ and 
$$C \forces_{\PP_2} \; h_\beta \in \dot D_\beta.$$
I.e.
$$(C,h_\beta) \forces_{\PP_2 * \QQ} \dot\Phi (\beta) \restrict \dot A
=\phi_\beta \restrict \dot A = \dot\Phi(\gamma) \restrict \dot A$$
for some $\gamma = \gamma(\beta) < \alpha$. Choose by fact (3) $Z'
\in [Z]^\kappa$ and $\tilde h \supseteq h$ ($\tilde h \in {\cal A}$ of course)
and $\gamma$ such that \par
\item{(1)} $\forall \beta_1 \ne \beta_2 \in Z' \; ( dom (h_{\beta_1})
\cap dom (h_{\beta_2}) = dom(\tilde h))$; 
$\forall \beta \in Z' \; ( \tilde h = h_\beta \restrict dom(\tilde h))$;
\par
\item{(2)} $\gamma(\beta) = \gamma$ for $\beta \in Z'$. \par
\noindent $C \leq W$ is a finite abelian subgroup.
So $[W : C_W(C)] < \omega$. As $\vert Z'\vert = \kappa$,
$\vert \{ \phi_\beta \restrict C_W(C) ; \; \beta \in Z' \} 
\vert = \kappa$ so that we can find $c \in C_W(C) \setminus C$ and
$\beta_1, \beta_2 \in Z'$ such that $ (c)\phi_{\beta_1} \ne 
(c)\phi_{\beta_2} $. Then
$$(\langle C,c \rangle, h_{\beta_1} \cup h_{\beta_2} )
\forces_{\PP_2 * \QQ} " \dot \Phi (\beta_1) \restrict \dot A = \phi_{\beta_1}
\restrict \dot A = \dot\Phi (\gamma) \restrict \dot A = \phi_{\beta_2}
\restrict \dot A=\dot\Phi (\beta_2) \restrict \dot A$$
$$ {\rm  and } \;
\phi_{\beta_1}\restrict \dot A \ne \phi_{\beta_2} \restrict \dot A,"$$
which is a contradiction. \par This ends the proof of the main claim
and shows that there is indeed an abelian $A \leq W$ such
that $[G:C_G(A)] = \omega_2$. Then either $[G:N_G(A)] =
\omega_2$ or we apply Theorem 3.2 -- as we did before in the proofs
of Theorems B and C -- to get $B \leq A$ such that $[G:N_G(B)]
= \omega_2$. This is the final contradiction. $\qed$
\bigskip
It should be clear that this proof also yields that in Mitchell's
model, both ${\cal Y}_{\omega_1} = {\cal Z}_{\omega_1}$ and
${\cal Y}_{\omega_2} = {\cal Z}_{\omega_2}$ for periodic $FC$-groups;
especially ${\cal Y} = {\cal Z}$ for groups of size $\leq \omega_2$.
\bigskip
{\sanse 5.8.} {\capit Theorem.} {\it The consistency of $ZFC + I$ implies
the consistency of $ZFC \; +$ the following statements. \par
\item{(i)} ${\cal Y}_{\omega_1} = {\cal Z}_{\omega_1}$ for 
periodic $FC$-groups -- and
any periodic $FC$-group $G$ with $\vert G / Z(G) \vert = \omega_1$ has an abelian subgroup
$A$ with $[G:N_G(A)] = \omega_1$; \par 
\item{(ii)} ${\cal Y}_{\omega_2} = {\cal Z}_{\omega_2}$ for periodic
finite-by-abelian groups -- and any periodic finite-by-abelian group
$G$ with $\vert G/Z(G) \vert = \omega_2$ has an abelian subgroup
$A$ with $[G:N_G(A)] = \omega_2$; \par
\item{(iii)} There is an $FC$-group $G$ with $\vert G/Z(G) \vert =
\omega_2$ but $[G:N_G(A)] \leq \omega_1$ for all abelian $A \leq G$. \par}
\smallskip
{\it Proof.} Let $\kappa$ be inaccessible in ${\cal V}$.
By [1, Theorem 8.8], there is a partial order $\PP$ such that
for ${\cal G}$ $\PP$-generic over ${\cal V}$, \smallskip
\centerline{${\cal V} [{\cal G}] \models MA + 2^\omega = \omega_2
+$ "there are no weak Kurepa trees".} \smallskip
\noindent (This is proved by a countable support iteration of length
$\kappa$ of partial orders which alternatively make $\omega_1$-trees
special and are ccc.) Now apply 4.4/4.5 (for (i)), 4.6 (for (iii)),
and 5.4 (for (ii)). $\qed$
\bigskip
As both $MA$ and {\it no weak Kurepa trees} follow from the proper
forcing axiom $PFA$ (by [1, $\S$ 8] -- see also [2, 7.10]),
(i) through (iii) in the Theorem hold if we assume $ZFC + PFA$.
\bigskip
{\sanse 5.9.} {\capit Proposition.} {\it For any cardinal $\kappa \geq
\omega_2$ it is consistent that there is an extraspecial $p$-group
of size $\kappa$ such that for all maximal abelian subgroups 
$A \leq G$, $[G:A] \leq \omega_1$.}
\smallskip
{\it Proof.} By the arguments of 5.2 and 5.3 it suffices to generically
add a Kurepa tree with $\kappa$ branches as follows (Folklore).
Let $\KK_\kappa := \{ p ; \; p$ is a function and $dom(p) = \alpha \times A
$, where $\alpha < \omega_1$ and $A \in [\kappa]^\omega$, and
$ran(p) \subseteq 2 \}$, ordered by $p \leq q$ iff
$\alpha(p) \geq \alpha(q)$ (where $\alpha(p)$ is the $\alpha$
of the definition of the p.o.), $A(p) \supseteq A(q)$,
$p \restrict (\alpha(q) \times A(q)) = q$, and  for all $\beta \in
A(p) - A(q)$ there is $\gamma \in A(q)$ such that $p \restrict (\alpha(q)
\times \{\beta\}) = p \restrict (\alpha(q) \times \{\gamma\})$.
$\KK_\kappa$ is (if we assume $CH$ in the ground model ${\cal V}$)
$\omega_2 - cc$ and $\omega_1$-closed, and so preserves cardinals.
Clearly this works. $\qed$
\Bigskip
\vfill\eject

\noindent{\dunhg $\S$ 6. Generalizations}
\Smallskip
{\sanse 6.1.} It was mentioned in the Introduction that (II) holds
for arbitrary ${\cal Z}_\kappa$-groups $G$ [4]. One might ask what goes
wrong in this general case (where $G$ is not required to be $FC$) if
we drop the ${\cal Z}_\kappa$-condition.
\smallskip
{\capit Proposition.} {\it For any cardinal $\kappa$ there is a group $G$
with $\vert G /Z(G) \vert = 2^\kappa$ and $[G:N_G(A) ] \leq \kappa$
for all abelian $A \leq G$.}
\smallskip
{\it Proof.} Let $A$ and $B$ be two elementary abelian $p$-groups
of size $\kappa$. Let $h: [\kappa]^2 \to \kappa$ be a bijection.
We define a factor system $\tau$ as follows:
$$\tau(a_\alpha, a_\beta) = \cases{0 &if $\alpha \geq \beta$, \cr
b_{h(\alpha,\beta)} &otherwise, \cr}$$
where the $a_\alpha$ ($b_\alpha$, resp.) ($\alpha < \kappa$)
generate $A$ ($B$, resp.); extend $\tau$ bilinearly to $A^2$.
Let $C:=E(\tau)$ be the extension. (Note that $C$ is the free object
on $\kappa$ generators in the variety of two-step nilpotent groups
of exponent $p$ ($p > 2$); and that its maximal abelian subgroups
are of the form $\langle B, a \rangle$, where $a \not\in B$.) 
Let $D$ be a group of order $p$. Let $G$
be the (abelian) subgroup of $Aut(C \oplus D)$ consisting of all automorphisms $\phi$
which fix $B \oplus D$ and satisfy
$$ \forall c \in C \oplus D \; \exists k \in p \; (c\phi = c + kd).$$
Clearly, $\vert G \vert = 2^\kappa$. Let $E$ be the semidirect extension
of $C \oplus D$ and $G$ (i.e. $E=(C \oplus D) \sd G$). We leave 
it to the reader to verify that $\vert E / Z(E)
\vert = 2^\kappa$ and $[E:C_E(A)] \leq \kappa$ for all abelian
$A \leq E$. $\qed$
\smallskip
{\it Note.} The proof is similar to (but easier than) the proof of
Theorem A (see 4.2 and 4.3). Unlike the latter it does not involve
any set-theoretic hypotheses. On the other hand, the group $E$ constructed
above is not $\kappa C$ but $\kappa^+ C$ (a group $G$ is $\kappa C$ iff
every $g \in G$ has less than $\kappa$ conjugates).
\bigskip
{\sanse 6.2.} We restricted our attention to $\kappa = \omega_1$ or $\omega_2$.
This is reasonable because the problem seems to be most interesting for
small cardinals. Also, the constructions in $\S 5$ (5.3 and 5.9) show
how to get consistency results concerning the existence of
pathological groups for larger cardinals (just use 
$\lambda$-Kurepa families instead of ($\omega_1$-)Kurepa families for 
the appropriate $\lambda$). Nevertheless we ignore whether 
the non-existence of such groups
is consistent for $\kappa \geq \omega_3$ (cf Theorem E).   
\vfill\eject

\noindent{\dunhg References}
\Smallskip

\item{1.} {\capit J. Baumgartner,} "Iterated forcing,"
Surveys in set theory (edited by A. R. D. Mathias), Cambridge University
Press, Cambridge, 1983, 1-59.
\smallskip
\item{2.} {\capit J. Baumgartner,} "Applications of the
proper forcing axiom," Handbook of set-theoretic topology, North-Holland,
Amsterdam, 1984, 913-959.
\smallskip
\item{3.} {\capit V. Faber, R. Laver and R. McKenzie,}
Coverings of groups by abelian subgroups, {\it Canad. J. Math.}
{\bf 30} (1978), 933-945.
\smallskip
\item{4.} {\capit V. Faber and M. J. Tomkinson,} On 
theorems of B. H. Neumann concerning $FC$-groups II, {\it Rocky Mountain
J. Math.} {\bf 13} (1983), 495-506.
\smallskip
\item{5.} {\capit L. Fuchs,} "Infinite Abelian Groups I,"
Academic Press, New York, 1970.
\smallskip
\item{6.} {\capit L. Fuchs,} "Infinite Abelian Groups II,"
Academic Press, New York, 1973.
\smallskip
\item{7.} {\capit T. Jech,} "Set theory," Academic Press,
San Diego, 1978.
\smallskip
\item{8.} {\capit T. Jech,} "Multiple forcing," 
Cambridge University Press, Cambridge, 1986.
\smallskip
\item{9.} {\capit K. Kunen,} "Set theory," North-Holland,
Amsterdam, 1980.
\smallskip
\item{10.} {\capit L. A. Kurdachenko,} Dvustupenno nil'potennye
FC-gruppy (Twostep nilpotent
FC-groups), {\it Ukrain. Mat. Zh.} 
{\bf 39} (1987), 329-335.
\smallskip
\item{11.} {\capit W. Mitchell,} Aronszajn trees and the
independence of the transfer property, {\it Ann Math. Logic}
{\bf 5} (1972), 21-46.
\smallskip
\item{12.} {\capit D. J. S. Robinson,} "A course on the theory 
of groups," Springer, New York Heidelberg Berlin, 1980. 
\smallskip
\item{13.} {\capit S. Shelah and J. Stepr\=ans,} 
Extraspecial $p$-groups, {\it Ann. Pure Appl. Logic}
{\bf 34} (1987), 87-97.
\smallskip
\item{14.} {\capit M. J. Tomkinson,} Extraspecial sections
of periodic $FC$-groups, {\it Compositio Math.} {\bf 31} (1975),
285-302.
\smallskip
\item{15.} {\capit M. J. Tomkinson,} On theorems of B. H. 
Neumann concerning $FC$-groups, {\it  Rocky Mountain J.
Math.} {\bf 11} (1981), 47-58.
\smallskip
\item{16.} {\capit M. J. Tomkinson,} "$FC$-groups,"
Pitman, London, 1984.
\vfill\eject\end